\documentstyle[a4,%11pt,%
twoside]{article}
\makeatletter
\renewcommand{\@oddhead}{An Interaction of An Oscillator with 
  An One-Dimensional Scalar Field. \hfill \thepage}
\renewcommand{\@evenhead}{\thepage \hfill S.A. Choro\v{s}avin }
\renewcommand{\@oddfoot}{}
\renewcommand{\@evenfoot}{}
\makeatother
{\vspace*{2ex}\par {#1}\quad\it }{\medskip\par }

\newenvironment{Thm}[2]%
{\par\addvspace{\bigskipamount}{\bf #1#2}\it }%
{\par\addvspace{\bigskipamount} }

{\par\hspace*{\fill}$\Box$\par\addvspace{\bigskipamount} }

{\par\hspace*{\fill}$\Box$ \par\addvspace{\bigskipamount} }

\newcommand{\intt}%
{
\makebox[3ex][c]
{\makebox[0pt]{$-$}\makebox[0pt]{$\displaystyle \int$}}
}

\author{S.A.~Choro\v{s}avin}
\title{ 
  An Interaction of An Oscillator with 
  An One-Dimensional Scalar Field.
 \bigskip\\
  Simple Exactly Solvable Models based on \\ 
  Finite Rank Perturbations Methods. \bigskip\\ 
  I: D'Alembert-Kirchhoff-like formulae }
\date{}
\begin{document}
%%%%%%%%%%%%%%%%%%%%%%%%%%%%%%%%%%%%%%%%%%%%%%%%%%%%%%%%%%%%%%%%%%%%%%%%%%%%%%%
%%%%%%%%%%%%%%%%%%%%%%%%%%%%%%%%%%%%%%%%%%%%%%%%%%%%%%%%%%%%%%%%%%%%%%%%%%%%%%%
%%%%%%%%%%%%%\hfill S. A. Skarst Choro\v{s}avin. \hfill \\
\maketitle 
\begin{abstract}
 This paper is an electronic application to my set of lectures, 
 subject:`Formal methods in solving differential equations and 
 constructing models of physical phenomena'. Addressed, mainly: 
 postgraduates and related readers. 
 Content: a very detailed discussion of the simple model of interaction based 
 on the equation array: 
% d^2 q(t)/dt^2 =-\Omega^2(q(t)-Q(t))+f_0(t) , 
% d^2 u(t,x)/dt^2=c^{2}d^2 u(t,x)/dx^2
%                 -4\gamma c\delta(x-x_0)(Q(t)-q(t))+f_1(t,x) , 
% Q(t) = u(t,x_0) .
\begin{eqnarray*}
 \frac{\partial^2 q(t)}{\partial t^2}
 &=&-\Omega^2\Big(q(t)-Q(t)\Big)+f_0(t) 
\\
 \frac{\partial^2 u(t,x)}{\partial t^2} 
 &=&
 c^{2}\frac{\partial^2 u(t,x)}{\partial x^2}
 -4\gamma c\delta(x-x_0)\Big(Q(t)-q(t)\Big)+f_1(t,x) 
\\
 Q(t)  
 &=& u(t,x_0)   
\\
\end{eqnarray*}
 Besides, less detailed discussion of related models.
 Central mathematical points: d'Alembert-Kirchhoff-like formulae. 
 Central physical points: phenomena of Radiation Reaction, Braking Radiation 
 and Resonance.
\end{abstract}

%\begin{center} 
%{\Large \bf 
%  Once More on Interaction of An Oscillator with 
%  An One-Dimensional Scalar Field.
% \bigskip\\
%  Simple Exactly Solvable Models based on \\ 
%  Finite Rank Perturbations Methods \bigskip\\ 
%  I: D'Alembert-Kirchhoff-like formulae }
%\end{center} 

\newpage 
\section*%
{ Introduction. } 
\subsection*%
{ A Harmonic Oscillator Coupled to an One-Dimensional Scalar Field. }

 In this paper I will discuss several models of 
 an one-dimensional harmonic oscillator 
 coupled to an one-dimensional scalar field. 
 Primarily I am interested in the model described by the equation array
\footnote{
$4\gamma c\rho = k\,, \Omega^2 = k/M$ ;  
 the constants of the model mean, e.g.:
$c=$ 
 propagating waves velocity, 
$\rho = $ 
 `a density' of the field, 
$k = $ 
 elasticity constant, 
$M = $ 
 mass of the particle.
 Of course, we assume 
$c>0$ . 
}
\begin{eqnarray*}
 \frac{\partial^2 q(t)}{\partial t^2}
 &=&-\Omega^2\Big(q(t)-Q(t)\Big)+f_0(t) 
\\
 \frac{\partial^2 u(t,x)}{\partial t^2} 
 &=&
 c^{2}\frac{\partial^2 u(t,x)}{\partial x^2}
 -4\gamma c\delta(x-x_0)\Big(Q(t)-q(t)\Big)+f_1(t,x) 
\\
 Q(t)  
 &=& u(t,x_0)   
\\
\end{eqnarray*}
 One can say that this is a model of a point interaction.
 The considered model is not standard, 
 if one means classical or quantum field models.
 In that field, a standard looks rather like this:
\begin{eqnarray*}
 \frac{\partial^2 q_0(t)}{\partial t^2}
 &=&-\Omega^2 q_0(t) + \gamma_1 Q_{\phi}(t) + f_0(t) 
\\
 \frac{\partial^2 \phi(t,x)}{\partial t^2} 
 &=&
 c^{2}\frac{\partial^2 \phi(t,x)}{\partial x^2}
 +4\gamma_2 c\delta(x-x_0)q_0(t)+f_1(t,x) 
\\
 Q_{\phi}(t)  
 &=& \phi(t,x_0)   
\\
\end{eqnarray*}
 After indicating {\bf d'Alembert-Kirchhoff-like formulae }
 for solutions to these systems I briefly compare them 
 and discuss 
% phenomenon of {\bf Braking Radiation}
% and one of {\bf Resonance}.
 phenomena of {\bf Radiation Reaction}, {\bf Braking Radiation} 
 and {\bf Resonance}.
 
 After that I will discuss an abstract analogue of the former system: 
 \footnote%
 { the proper modifications for the latter system will be evident }
\begin{eqnarray*}
 \ddot q
 &=&-\Omega^2(q-Q)+f_0(t) 
\\
 \ddot u 
 &=&
 Bu -4\gamma_c\Big(\delta_{\alpha,t,x_0}\Big)\cdot\Big(Q-q\Big)+f_1(t) 
\\ 
 Q  
 &=& Q(t) = <l|u(t)> 
\\
\end{eqnarray*}
 where I use a P.A.M. Dirac's ``bra-ket'' syntax and suppose that 
$q$ and $Q$
 are usual (one-dimensional) functions of 
$t$ : 
$$
  q=q(t) \,,\quad Q=Q(t) \,, 
$$ 
$B$ 
 is an abstract linear operator,
$l$
 is a linear functional,
$\{u(t)\}_{t}$ 
 and 
$\{\delta_{\alpha,t,x_0}\}_{t}$ 
 are families of abstract elements; 
 of course the type of 
$\delta_{\alpha,t,x_0}$ 
 must be the same as one of 
$u(t)$.

 One can rewrite the above equations as follows: 
$$
\left(\begin{array}{cc}
 \ddot q \\
 \ddot u
\end{array}\right)
 =
\left(\begin{array}{cc}
 -\Omega^2         
       &  \Omega^2<l|_2 \\
 4\gamma_c\delta_{\alpha,t,x_0}<1|_1 
       &  B -4\gamma_c\delta_{\alpha,t,x_0}<l|_2
\end{array}\right)
\left(\begin{array}{c}
 q \\
 u
\end{array}\right)
 + 
\left(\begin{array}{cc}
 f_0(t) \\
 f_1(t)
\end{array}\right)
$$
$$
\left(\begin{array}{cc}
 \ddot q_0 \\
 \ddot \phi 
\end{array}\right)
 =
\left(\begin{array}{cc}
 -\Omega^2         
       &  \gamma_1<l|_2 \\
 4\gamma_{2,c}\delta_{\alpha,t,x_0}<1|_1 
       &  B 
\end{array}\right)
\left(\begin{array}{c}
 q_0 \\
 \phi 
\end{array}\right)
 + 
\left(\begin{array}{cc}
 f_0(t) \\
 f_1(t)
\end{array}\right)
$$
 where the subscribts 
$1$ and $2$ in 
$<\cdots|_1$ and $<\cdots|_2$ 
 mean that arguments of 
$<\cdots|_1$ 
 are elements of the first component of the vector 
$\left(\begin{array}{c}
 q \\
 u
\end{array}\right)$
 resp.
$\left(\begin{array}{c}
 q_0 \\
 \phi
\end{array}\right)$
 and arguments of 
$<\cdots|_2$ 
 are elements of the second component of the suitable vector.
 Of course, 

$$ 
 <1|q> = <1|_1 q> = q\,, \quad <1|q_0> = <1|_1 q_0> = q_0 
$$

 Normally I will suppose that 
$\delta_{\alpha,t,x_0}$ 
 does not depend on 
$t$, 
 i.e. is constant in 
$t$. 
 In that case I will write  
$\delta_{\alpha,x_0}$ 
 instead of 
$\delta_{\alpha,t,x_0}$.

 In the next paper I will 
 suppose that 
$\delta_{\alpha,t,x_0}$ 
 is constant in 
$t$. 

 The subject will primarily be 
 {\bf resolvents formulae}
 i.e. 
 the formulae that resolve the equation array:  
$$
 z
\left(\begin{array}{cc}
  q \\
  u
\end{array}\right)
 -
\left(\begin{array}{cc}
 -\Omega^2         
       &  \Omega^2<l|_2 \\
 4\gamma_c\delta_{\alpha,x_0}<1|_1 
       &  B -4\gamma_c\delta_{\alpha,x_0}<l|_2
\end{array}\right)
\left(\begin{array}{c}
 q \\
 u
\end{array}\right)
 = 
\left(\begin{array}{cc}
 w_1 \\
 w_2
\end{array}\right)
$$
 resp.
$$
 z
\left(\begin{array}{cc}
  q_0 \\
  \phi 
\end{array}\right)
 -
\left(\begin{array}{cc}
 -\Omega^2         
       &  \gamma_1<l|_2 \\
 4\gamma_{2,c}\delta_{\alpha,x_0}<1|_1 
       &  B 
\end{array}\right)
\left(\begin{array}{c}
 q_0 \\
 \phi 
\end{array}\right)
 = 
\left(\begin{array}{cc}
 w_1 \\
 w_2
\end{array}\right)
$$

 Of course, here all quantities,
$q, u, w_1, w_2,\delta_{\alpha,x_0}, l $, 
 etc. are supposed to be constant in  
$t$. 

 I will discuss 
 Donoghue-Friedrichs-like solutions to the system
 \footnote{ i.e., I will discuss resolvents formulae of the system } 
 and resp.
 the phenomenon of {\bf Resonance} and notion of the {\bf Second Sheet}.

\newpage\section%*%
{ Models of a Point Interaction of an only Oscillator with 
  an only one-dimensional Scalar Field }

 In this section we fix measure units 
% so that 
%$c >0$ 
 and let
$x$ 
 be dimensionless position parameter, i.e., 
$$
 \mbox{ physical position coordinate } 
  =  [\mbox{ length unit }]\times x +const \,. 
$$
 Otherwise a confusion can ocurr, in relating to the definition 
$$
 \int_{-\infty}^{\infty}\delta(x-x_0)f(x)dx=f(x_0) \,.
$$
 We assume the standard foramalism, where 
$$ 
 \delta(x-x_0) = \frac{\partial 1_{+}(x-x_0)}{\partial x} 
$$
 and where 
$1_{+}$
 stands for a unit step function (Heaviside function):
$$ 
 1_{+}(\xi) = \Bigl\{ 
 \begin{array}{ccc}
 1&,&\mbox{ if } \xi > 0 \,,
 \\
 0&,&\mbox{ if } \xi < 0 \,,
 \\ 
 \end{array}
$$           
%\newpage
\subsection
{ First Model. D'Alembert-Kirchhoff-like formulae }

 Recall that standard D'Alembert-Kirchhoff formulae read: 
 if 
$$
 \frac{\partial^2 u}{\partial t^2}
  = c^2\frac{\partial^2 u}{\partial x^2} + f
 \,,\quad
 u=u(t,x)
 \,,\quad
 f=f(t,x)
 \,,\quad
% (t\geq s)
           \eqno(*)
$$
 and given initial data, 
$u(s,\cdot )$ 
 and 
$\frac{\partial u(t,\xi)}{\partial t})\Big|_{t=s}$,
 then 
\begin{eqnarray*}
 u=u(t,x)&=&
 \frac{1}{2c}
 \int_s^t 
 \Big(
 \widetilde f(\tau,x+c(t-\tau))-\widetilde f(\tau,x-c(t-\tau))\Big)
 d\tau 
\\&&{}
 + u_{0}(t,x)
\end{eqnarray*}
%%%
\begin{eqnarray*}
 u_{0}(t,x)
 &=& 
 c_{+}(x+ct) + c_{-}(x-ct)
 \\&=&
 \frac12 \Big(u(s,x+ct) + u(s,x-ct)\Big)
 \\&&{}
   + \frac{1}{2c}\Big(
  \widetilde{\dot u}(s,x+ct) - \widetilde{\dot u}(s,x-ct)\Big)
 \\&&{}
\end{eqnarray*}
 and where 
$\widetilde f \,,\, \widetilde{\dot u}$ 
 stand for any functions defined by 
$$
  \frac{\partial\widetilde f(t,x)}{\partial x} = f(t,x)
\,,\quad 
  \frac{\partial\widetilde{\dot u}(s,\xi)}{\partial \xi}
    =\Bigl(\frac{\partial u(t,\xi)}{\partial t}\Bigr)\Big|_{t=s}
\,.
$$
 Note that
$$
 \widetilde f(\tau,x+c(t-\tau))-\widetilde f(\tau,x-c(t-\tau))
\,,\quad 
\widetilde{\dot u}(s,x+ct) - \widetilde{\dot u}(s,x-ct)
$$
 do not depend on what the primitives are which one has chosen!!! 
 Moreover, we need only  
$\widetilde f \,,\, \widetilde{\dot u}|_{t=s} $
 and not 
$f \,,\, \dot u|_{t=s}$
 themselves!

 Recall that standard relation for an one-dimensional harmonic oscillator
 under an external force 
$f_0(t)$ 
 is this: 
$$
 \frac{\partial^2 q(t)}{\partial t^2}
  =-\Omega^2 q(t)+f_0(t) 
$$

 In this section I will firstly discuss a system described by  

\begin{eqnarray*}
 \frac{\partial^2 q(t)}{\partial t^2}
 &=&-\Omega^2\Big(q(t)-Q(t)\Big)+f_0(t) 
\\
 \frac{\partial^2 u(t,x)}{\partial t^2} 
 &=&
 c^2\frac{\partial^2 u(t,x)}{\partial x^2}
 -4\gamma c\delta(x-x_0)\Big(Q(t)-q(t)\Big)+f_1(t,x) 
\\
 Q(t)  
 &=& u(t,x_0)   
\\
\end{eqnarray*}
 It means in patricular that I will take  
$$
  f=-4\gamma c\delta(x-x_0)\Big(Q(t)-q(t)\Big)+f_1(t,x) 
$$
 For such an  
$f$  
 I conclude that 
$$
 \tilde f=-4\gamma c1_{+}(x-x_0)\Big(Q(t)-q(t)\Big)+\tilde f_1(t,x), 
$$
%
%%%%%     \newpage 
%
 and then I infer that 
\begin{eqnarray*}
 u(t,x)
 &=&
 -{2\gamma}
 \int_s^t \Big( 1_{+}(x+c(t-\tau)-x_0)-1_{+}(x-c(t-\tau)-x_0)\Big)
 \Big(Q(\tau)-q(\tau)\Big)d\tau
 \\&&{} + \frac{1}{2c}\int_s^t 
   \Big(\tilde f_1(\tau,x+c(t-\tau))-\tilde f_1(\tau,x-c(t-\tau))\Big) d\tau 
% \\&&{} 
 +  u_{0}(t,s,x)
%c_{+}(x+ct) + c_{-}(x-ct)
\\
 u(t,x)
 &=&
 -{2\gamma}
 \int_s^t \Big( 1_{+}(x+c(t-\tau)-x_0)-1_{+}(x-c(t-\tau)-x_0)\Big)
 \Big(u(\tau,x_0)-q(\tau)\Big)d\tau
 \\&&{} + \frac{1}{2c}\int_s^t 
   \Big(\tilde f_1(\tau,x+c(t-\tau))-\tilde f_1(\tau,x-c(t-\tau))\Big) d\tau 
% \\&&{}
 +  u_{0}(t,s,x)
%c_{+}(x+ct) + c_{-}(x-ct)
\end{eqnarray*}
 Denote now, to be more concise, 
$$ 
 u_{01}(t,s,x)
 := 
  \frac{1}{2c}\int_s^t 
  \Big(\tilde f_1(\tau,x+c(t-\tau))-\tilde f_1(\tau,x-c(t-\tau))\Big) d\tau 
  +  u_{0}(t,s,x)
%c_{+}(x+ct) + c_{-}(x-ct)
$$ 
 and then rewrite the recent relation %equation for 
%$u(t,x)$ 
 as following: 
\begin{eqnarray*}
 u(t,x)
 &=&
 -{2\gamma}
 \int_s^t \Big( 1_{+}(x+c(t-\tau)-x_0)-1_{+}(x-c(t-\tau)-x_0)\Big)
 \Big(u(\tau,x_0)-q(\tau)\Big)d\tau
 \\&&{} + u_{01}(t,s,x) 
\end{eqnarray*}
 I have now seen: given
$q$
 and 
$u_{01}$,
 then, 
\bigskip 

\fbox{\bf 
 in order to obtain 
$ u(t,x) $
 I need to obtain ONLY 
$ u(t,x_0) $ } 
\bigskip 
\\
 After this observation use the last formula for 
$u(t,x)$
  and then obtain
\begin{eqnarray*}
 u(t,x_0)
  &=&
 -{2\gamma}
 \int_s^t 
 \Big(1_{+}(c(t-\tau)) - 1_{+}(c(\tau-t))\Big)
 \Big(u(\tau,x_0)-q(\tau)\Big)d\tau 
 \\&&{} + u_{01}(t,s,x_0) 
\end{eqnarray*}
 Therefore,  because 
$$
   1_{+}(c(t-\tau))=1\,,\, 1_{+}(c(\tau-t)) = 0 \mbox{ for } t>\tau 
$$
 one can obtain

\begin{eqnarray*}
 u(t,x_0)
 &=&
 -{2\gamma}
 \int_s^t 
 \Big(u(\tau,x_0)-q(\tau)\Big)d\tau 
 + u_{01}(t,s,x_0) 
\end{eqnarray*}
 i.e.,
\begin{eqnarray*}
 Q(t)
 &=&
 -{2\gamma}
 \int_s^t 
 \Big(Q(\tau)-q(\tau)\Big)d\tau 
 + u_{01}(t,s,x_0) 
\end{eqnarray*}
 Write it as 

\begin{eqnarray*}
 Q(t) + {2\gamma} \int_s^t Q(\tau) d\tau
 &=&
 {2\gamma} \int_s^t q(\tau) d\tau
 + u_{01}(t,s,x_0) 
\end{eqnarray*}
 Now recall that 
$$
 \frac{\partial^2 q(t)}{\partial t^2} 
 = -\Omega^2 \Big(q(t)-Q(t)\Big) + f_0(t)  
$$
 and for the moment denote
$$
 Q_{0}(t) := u_{0}(t,s,x_0) 
           = c_{+}(x_0+ct) + c_{-}(x_0-ct)
$$
$$
 Q_{01}(t) 
 := u_{01}(t,s,x_0)
  =  \frac{1}{2c}\int_s^t 
     \Big(
     \tilde f_1(\tau,x_0+c(t-\tau))-\tilde f_1(\tau,x_0-c(t-\tau))
     \Big) d\tau 
     + Q_{0}(t)
$$
 Then obtain
\begin{eqnarray*} 
 Q(t) + {2\gamma} \int_s^t Q(\tau) d\tau
 &=&
 {2\gamma} \int_s^t q(\tau) d\tau
 + Q_{01}(t) 
 \\
  \frac{\partial^2 q}{\partial t^2} 
 &=& -\Omega^2(q-Q)  + f_0(t) 
\end{eqnarray*}

 We have now obtained an equation array for 
$Q$ and $q$. The next step is to find an {\bf insulated} equation for 
$Q$ 
 and one for
$q$. 
 We begin to search for the equation for
$q$. 
 For this purpose, we apply, at first, the operator 
$$
  \Big( I + {2\gamma} \int_s^t \cdot d\tau  \Big)
$$
 to the latter equation, i.e. to the equation 
\begin{eqnarray*} 
  \frac{\partial^2 q}{\partial t^2} 
 &=& -\Omega^2(q-Q)  + f_0(t) 
\end{eqnarray*}
 Then we infer 
\begin{eqnarray*}
  \frac{\partial^2 q}{\partial t^2} 
 +{2\gamma} \int_s^t \frac{\partial^2 q(\tau)}{\partial \tau^2}d\tau  
 &=&
 -\Omega^2 \Big( q +{2\gamma} \int_s^t q(\tau) d\tau \Big)
\\&&{} 
 + \Omega^2 \Big( Q +{2\gamma} \int_s^t Q(\tau) d\tau \Big)
\\&&{}
 + \Big( f_0 +{2\gamma} \int_s^t f_0(\tau) d\tau \Big)
\\
  \frac{\partial^2 q}{\partial t^2} 
 +{2\gamma} \int_s^t \frac{\partial^2 q(\tau)}{\partial \tau^2}d\tau  
 &=&
 -\Omega^2 \Big( q +{2\gamma} \int_s^t q(\tau) d\tau \Big)
\\&&{}
 + \Omega^2 \Big(  {2\gamma} \int_s^t q(\tau) d\tau
 + Q_{01}(t) 
 \Big)
\\&&{}
 + \Big( f_0 +{2\gamma} \int_s^t f_0(\tau) d\tau \Big)
\\
  \frac{\partial^2 q}{\partial t^2} 
 +{2\gamma} \int_s^t \frac{\partial^2 q(\tau)}{\partial \tau^2}d\tau  
 &=&
 -\Omega^2
 \Big( q + Q_{01}(t) \Big)
 +
 \Big( f_0 +{2\gamma} \int_s^t f_0(\tau) d\tau \Big)
\\
\end{eqnarray*}
 The most recent equation for
$q$
 rewrite as
$$
 \frac{\partial^2 q}{\partial t^2}  
  = 
   -{2\gamma} 
    \int_s^t  \frac{\partial^2 q(\tau)}{\partial \tau^2}d\tau  
   -\Omega^2 q 
   + \Omega^2 Q_{01}(t)
   + f_0(t) 
   +{2\gamma} \int_s^t f_0(\tau) d\tau 
$$ 
 We can stop at this equation, or, observing that 
$$ 
 \int_s^t  \frac{\partial^2 q(\tau)}{\partial \tau^2}d\tau  
 = \frac{\partial q(t)}{\partial t}
   -\frac{\partial q(t)}{\partial t}\Big|_{t=s}
 \,, 
$$
 we can stop at that:

\medskip 

\hspace{\fill}
\fbox{
\parbox{\textwidth}{
$$
 \frac{\partial^2 q(t)}{\partial t^2} 
  =-{2\gamma}\frac{\partial q(t)}{\partial t} 
   -\Omega^2 q(t) 
   +{2\gamma}\frac{\partial q(t)}{\partial t}\Big|_{t=s} 
   + \Omega^2 Q_{01}(t) 
   + f_0(t) 
   +{2\gamma} \int_s^t f_0(\tau) d\tau 
$$ 
}%%% end of parbox
}%%% end of fbox
\hspace{\fill}
\medskip 

\noindent 
  Some people prefer to write such an equation as following: 

\medskip 

$$
\fbox{
\parbox{\textwidth}{
$$
 \Bigl(\frac{\partial^2 }{\partial t^2} 
  +{2\gamma}\frac{\partial }{\partial t} 
  +\Omega^2 \Bigr)q(t) 
  =
   {2\gamma}\frac{\partial q(t)}{\partial t}\Big|_{t=s} 
   + \Omega^2 Q_{01}(t) 
   + f_0(t) 
   +{2\gamma} \int_s^t f_0(\tau) d\tau 
$$ 
}%%% end of parbox
}%%% end of fbox
$$

\bigskip 

\noindent  
 Remark, whatever equation we take, we see an interesting detail:
 the equation is {\bf not ordinary } differential equation,
 if we follow standard terminology. The case is  
 because of the term 
$$
  {2\gamma}\frac{\partial q(t)}{\partial t}\Big|_{t=s} 
$$
 In the next subsections, we will return to this factor.  
 However, the machinery of the {\bf ordinary} differential equations 
 does here quite for.
 For the time being, we turn to describing 
$Q(t)$ 
 and 
$u(t,x)$.
\medskip 
 
 After  
$q(t)$
 is found,
 we can determine 
$Q(t)$ ,
 at least formally, by solving 
\begin{eqnarray*}
 Q(t)
 &=&
 -{2\gamma}
 \int_s^t 
 \Big(Q(\tau)-q(\tau)\Big) d\tau 
 + Q_{01}(t)
 \\
\end{eqnarray*}
 or 
\begin{eqnarray*}
 \frac{\partial Q(t)}{\partial t}
 &=&
 -{2\gamma}
 \Big(Q(\tau)-q(\tau)\Big)
 + \frac{\partial Q_{01}(t)}{\partial t}
 \,,\quad Q(s) = Q_{01}(s)
 \\
\end{eqnarray*}
 
 Thus we have already reduced our model, a model of an oscillator coupled the 
 a scalar field, to a pair of linear `ordinary' differential equations. 
 Nevertheless we want to continue to analyse the matter 
 and we now go seeking another relationships, 
 which would simplify calculations of 
$q, Q$
 and 
$u$.

%%%% \newpage  

%
 At first, we will obtain another insulated equation for 
$Q$,
 differently and in a different form. 

 We have 
\begin{eqnarray*} 
 Q(t) + {2\gamma} \int_s^t Q(\tau) d\tau
 &=&
 {2\gamma} \int_s^t q(\tau) d\tau
 + Q_{01}(t) 
 \\
  \frac{\partial^2 q}{\partial t^2} + \Omega^2 q
 &=& \Omega^2 Q  + f_0(t) 
\end{eqnarray*}
 i.e.,
\begin{eqnarray*} 
 \Big(Q(t) - Q_{01}(t)\Big) 
 + {2\gamma} \int_s^t \Big(Q(\tau)-Q_{01}(\tau)\Big)d\tau 
 &=&
 {2\gamma} \int_s^t\Big(q(\tau)-Q_{01}(\tau)\Big) d\tau 
 \\
  \Big(\frac{\partial^2 }{\partial t^2} + \Omega^2 \Big)q 
 &=& \Omega^2 Q  + f_0(t) 
 \\
\end{eqnarray*}
 Then we infer 
\begin{eqnarray*} 
\makebox[30ex][l]{$\displaystyle 
 \Big(\frac{\partial^2 }{\partial t^2} + \Omega^2 \Big) 
 \Big(\Big(Q(t)-Q_{01}(t)\Big)
  + {2\gamma} \int_s^t \Big(Q(\tau)-Q_{01}(\tau)\Big) d\tau \Big) 
  $}%% end of the makebox 
\\&=&
 \Big(\frac{\partial^2 }{\partial t^2} + \Omega^2 \Big)
 \Big({2\gamma} \int_s^t \Big(q(\tau)-Q_{01}(\tau)\Big) d\tau \Big)
 \\
\end{eqnarray*}
 and 
\begin{eqnarray*} 
\makebox[15ex][l]{$\displaystyle 
 \Big(\frac{\partial^2 }{\partial t^2} + \Omega^2 \Big) 
 \Big(\Big(Q(t)-Q_{01}(t)\Big)
  + {2\gamma} \int_s^t \Big(Q(\tau)-Q_{01}(\tau)\Big) d\tau \Big) 
  $}%% end of the makebox 
\\&=&
 \Big(\frac{\partial^2 }{\partial t^2} + \Omega^2 \Big)
 \Big({2\gamma} \int_s^t q(\tau) d\tau \Big)
 - 
 \Big(\frac{\partial^2 }{\partial t^2} + \Omega^2 \Big)
 \Big({2\gamma} \int_s^t Q_{01}(\tau) d\tau \Big)
 \\
\end{eqnarray*}
 Note 
$$
 \frac{\partial^2 }{\partial t^2} \int_s^t q(\tau) d\tau 
 =
 \frac{\partial }{\partial t} q(t)  
 =
 \int_s^t\frac{\partial^2 }{\partial \tau^2} q(\tau) d\tau 
 +
 \frac{\partial }{\partial t} q(t) \Big|_{t=s} 
$$
 Then 
\begin{eqnarray*} 
\makebox[3ex][l]{$\displaystyle 
 \Big(\frac{\partial^2 }{\partial t^2} + \Omega^2 \Big) 
 \Big(\Big(Q(t)-Q_{01}(t)\Big)
  + {2\gamma} \int_s^t \Big(Q(\tau)-Q_{01}(\tau)\Big) d\tau \Big) 
  $}%% end of the makebox 
\\&=&
 {2\gamma} \int_s^t  
 \Big(\frac{\partial^2 }{\partial \tau^2} + \Omega^2 \Big)q(\tau) d\tau 
 +
 {2\gamma}\frac{\partial }{\partial t} q(t) \Big|_{t=s} 
 - 
 \Big(\frac{\partial^2 }{\partial t^2} + \Omega^2 \Big)
 \Big({2\gamma} \int_s^t Q_{01}(\tau) d\tau \Big)
 \\
\end{eqnarray*}
 Then 
\begin{eqnarray*} 
\makebox[2ex][l]{$\displaystyle 
 \Big(\frac{\partial^2 }{\partial t^2} + \Omega^2 \Big) 
 \Big(\Big(Q(t)-Q_{01}(t)\Big)
  + {2\gamma} \int_s^t \Big(Q(\tau)-Q_{01}(\tau)\Big) d\tau \Big) 
  $}%% end of the makebox 
\\&=&
 {2\gamma} \int_s^t  
 \Big( \Omega^2 Q(\tau)  + f_0(\tau) \Big) d\tau 
 +
 {2\gamma}\frac{\partial }{\partial t} q(t) \Big|_{t=s} 
 - 
 \Big(\frac{\partial^2 }{\partial t^2} + \Omega^2 \Big)
 \Big({2\gamma} \int_s^t Q_{01}(\tau) d\tau \Big)
 \\
\end{eqnarray*}
 Then 
\begin{eqnarray*} 
\makebox[3ex][l]{$\displaystyle 
 \Big(\frac{\partial^2 }{\partial t^2} + \Omega^2 \Big) 
 \Big(\Big(Q(t)-Q_{01}(t)\Big)
  + {2\gamma} \int_s^t \Big(Q(\tau)-Q_{01}(\tau)\Big) d\tau \Big) 
  $}%% end of the makebox 
\\&=&
 {2\gamma} \int_s^t  
 \Big(\Omega^2 Q(\tau)-\Omega^2 Q_{01}(\tau)\Big) d\tau 
 +
 {2\gamma}\frac{\partial }{\partial t} q(t) \Big|_{t=s} 
 - 
 \frac{\partial^2 }{\partial t^2}
 \Big({2\gamma} \int_s^t Q_{01}(\tau) d\tau \Big)
 \\&&{} \qquad \qquad
 + {2\gamma} \int_s^t f_0(\tau) d\tau
 \\
\end{eqnarray*}
 Then 
\begin{eqnarray*} 
\makebox[3ex][l]{$\displaystyle 
 \Big(\frac{\partial^2 }{\partial t^2} + \Omega^2 \Big) 
 \Big(Q(t)-Q_{01}(t)\Big)
  + 
 \Big(\frac{\partial^2 }{\partial t^2} + \Omega^2 \Big) 
 {2\gamma} \int_s^t \Big(Q(\tau)-Q_{01}(\tau)\Big) d\tau 
  $}%% end of the makebox 
\\&=&
 {2\gamma} \int_s^t  
 \Big(\Omega^2 Q(\tau)-\Omega^2 Q_{01}(\tau)\Big) d\tau 
 +
 {2\gamma}\frac{\partial }{\partial t} q(t) \Big|_{t=s} 
 - 
 \frac{\partial^2 }{\partial t^2}
 \Big({2\gamma} \int_s^t Q_{01}(\tau) d\tau \Big)
 \\&&{} \qquad \qquad
 + {2\gamma} \int_s^t f_0(\tau) d\tau
 \\
\end{eqnarray*}
 Then 
\begin{eqnarray*} 
\makebox[20ex][l]{$\displaystyle 
 \Big(\frac{\partial^2 }{\partial t^2} + \Omega^2 \Big) 
 \Big(\Big(Q(t)-Q_{01}(t)\Big)
  + 
 \frac{\partial^2 }{\partial t^2}
 {2\gamma} \int_s^t \Big(Q(\tau)-Q_{01}(\tau)\Big) d\tau \Big) 
  $}%% end of the makebox 
\\&=&
 {2\gamma}\frac{\partial }{\partial t} q(t) \Big|_{t=s} 
 - 
 \frac{\partial^2 }{\partial t^2}
 \Big({2\gamma} \int_s^t Q_{01}(\tau) d\tau \Big)
 + {2\gamma} \int_s^t f_0(\tau) d\tau
 \\
\end{eqnarray*}
 Then, finally,  
\medskip 

\hspace{\fill}
\fbox{
\hspace{\fill}
\parbox{1.01\textwidth}{
\begin{eqnarray*} 
 \Big(\frac{\partial^2 }{\partial t^2}
 +{2\gamma}\frac{\partial }{\partial t}
 + \Omega^2 \Big) 
 \Big(Q(t)-Q_{01}(t)\Big)
 &=&
 {2\gamma}\frac{\partial q(t)}{\partial t}\Big|_{t=s} 
 - 
 {2\gamma}\frac{\partial  Q_{01}(t)}{\partial t}
 + {2\gamma} \int_s^t f_0(\tau) d\tau
 \\
\end{eqnarray*}
}%%% end of parbox
\hspace{\fill}
}%%% end of fbox
\hspace{\fill}

\medskip 

\noindent 
 On the surface, this equation appears to be a second order {\bf ordinary } 
 differential equation. It is not exactly the case. 
 We may not arbitrary take the initial data
 for
$Q(t)$.
 The proper ones are these: 
$$
  \Big(Q(s)-Q_{01}(s)\Big) = 0\,,\,
  \frac{\partial\Big(Q(t)-Q_{01}(t)\Big)}{\partial t}\Big|_{t=s}
 = {2\gamma}\Big(q(s)-Q_{01}(s)\Big) 
$$%
%

%
%%%%%%%%%% \newpage    
%

 Thus, 
 we have already obtained two simple `ordinary' differential equations for 
$$
  q(t) \,,\quad Q(t)  
$$ 
 Now let us analyse the expression 
\begin{eqnarray*}
 u=u(t,x) 
 &=&
 -{2\gamma}
 \int_s^t \Big( 1_{+}(x+c(t-\tau)-x_0)-1_{+}(x-c(t-\tau)-x_0)\Big)
 \Big(Q(\tau)-q(\tau)\Big)d\tau
 \\&&{} + u_{01}(t,s,x) 
\end{eqnarray*}

 We have: if 
$
 c\tau \not= ct+(x-x_0) 
$
 and 
$
 c\tau \not= ct-(x-x_0) 
$
 then 
\begin{eqnarray*}
\makebox[5ex][l]{$\displaystyle 
 1_{+}(x+c(t-\tau)-x_0)-1_{+}(x-c(t-\tau)-x_0)
  $}%% end of the makebox 
\\ &=&
 \Big\{ 
 \begin{array}{rcl}
  1 & , & c\tau < ct+(x-x_0) \\ 
  0 & , & ct+(x-x_0) < c\tau \\ 
 \end{array}
 - 
 \Big\{ 
 \begin{array}{rcl}
  1 & , & ct-(x-x_0) < \tau \\ 
  0 & , & c\tau < ct-(x-x_0) \\ 
 \end{array}
\\ &=&
 \Bigg\{ 
 \begin{array}{rcl}
  1 & , & c\tau < ct-|x-x_0| \\ 
  0 & , & ct-|x-x_0| < c\tau < ct+|x-x_0| \\ 
 -1 & , & ct+|x-x_0| < c\tau \\ 
 \end{array}
\\ &=&
 \Bigg\{ 
 \begin{array}{rcl}
  1 & , & \tau < t-|x-x_0|/c \\ 
  0 & , & t-|x-x_0|/c < \tau < t+|x-x_0|/c \\ 
 -1 & , & t+|x-x_0|/c < \tau \\ 
 \end{array}
\end{eqnarray*}
 Hence 
\parbox{\textwidth}{
\begin{eqnarray*}
 u(t,x) 
 &=&
 \left\{ 
 \begin{array}{ccl}
 \displaystyle
 -{2\gamma}
 \int_s^{t-|x-x_0|/c} 
 \Big(Q(\tau)-q(\tau)\Big)d\tau &,& \mbox{ if } s \leq t-|x-x_0|/c \\ 
                 0              &,& \mbox{ if } t-|x-x_0|/c < s \leq t \\ 
                 0              &,& \mbox{ if } t \leq s < t+|x-x_0|/c \\ 
 \displaystyle
 +{2\gamma}
 \int_s^{t+|x-x_0|/c} 
 \Big(Q(\tau)-q(\tau)\Big)d\tau &,& \mbox{ if } t+|x-x_0|/c \leq s \\ 
 \end{array}
 \right\} 
 \\[\medskipamount]\\
 &&{}\qquad + u_{01}(t,s,x) 
 \\[\smallskipamount] \\ 
\end{eqnarray*}
}%%%end of parbox 
%
% Here we ought to say, 
% that when we deduced the latter formula, we had used the fact that 

\noindent 
 Besides, remember that 

\begin{eqnarray*}
 Q(t) + {2\gamma} \int_s^t Q(\tau) d\tau
 &=&
 {2\gamma} \int_s^t q(\tau) d\tau
 + Q_{01}(t) 
\end{eqnarray*}
  If we take this factor into account, then we deduce, finally, 

\parbox{\textwidth}{
\begin{eqnarray*}
 u(t,x) 
 &=&
 \left\{ 
 \begin{array}{ccl}
 \displaystyle
 Q(t-|x-x_0|/c)-Q_{01}(t-|x-x_0|/c) 
                                  &,& \mbox{ if } s \leq t-|x-x_0|/c \\ 
                 0                &,& \mbox{ if } t-|x-x_0|/c < s \leq t \\ 
                 0                &,& \mbox{ if } t \leq s < t+|x-x_0|/c \\ 
 -Q(t+|x-x_0|/c)+Q_{01}(t+|x-x_0|/c)
                                  &,& \mbox{ if } t+|x-x_0|/c \leq s \\ 
 \end{array}
 \right\} 
 \\[\medskipamount]\\
 &&{}\qquad + u_{01}(t,s,x) 
\end{eqnarray*}
}%%%end of parbox 

 Now substitute anywhere 
$ u_{0}(t,s,x_0) $ 
 for 
$ Q_{0}(t) $, 
 and 
$ u_{01}(t,s,x_0) $ 
 for
$ Q_{01}(t) $, 
 and 
 resume. 

\newpage\subsubsection
  {\bf Putting It Together } 

 We have discussed the system 
\begin{eqnarray*}
 \frac{\partial^2 q(t)}{\partial t^2}
 &=&-\Omega^2\Big(q(t)-Q(t)\Big)+f_0(t) 
\\
 \frac{\partial^2 u(t,x)}{\partial t^2} 
 &=&
 c^2\frac{\partial^2 u(t,x)}{\partial x^2}
 -4\gamma c\delta(x-x_0)\Big(Q(t)-q(t)\Big)+f_1(t,x) 
\\
 Q(t)  
 &=& u(t,x_0)   
\\
\end{eqnarray*}
 and concluded that:

\begin{eqnarray*}
 \Bigl(\frac{\partial^2 }{\partial t^2} 
  +{2\gamma}\frac{\partial }{\partial t} 
  +\Omega^2 \Bigr)q(t) 
 &=&
   {2\gamma}\frac{\partial q(t)}{\partial t}\Big|_{t=s} 
   + \Omega^2 u_{01}(t,s,x_0)
   + f_0(t) 
   +{2\gamma} \int_s^t f_0(\tau) d\tau 
\end{eqnarray*}
\begin{eqnarray*}
%\\ 
 Q(t)
 &=&
 -{2\gamma}
 \int_s^t 
 \Big(Q(\tau)-q(\tau)\Big) d\tau 
 + u_{01}(t,s,x_0)
 \\
\end{eqnarray*}
%
% The latter relation, one can rewrite it as the differential one:
%\begin{eqnarray*}
% \frac{\partial Q(t)}{\partial t}
% &=&
% -{2\gamma}
% \Big(Q(\tau)-q(\tau)\Big)
% + \frac{\partial u_{01}(t,s,x_0)}{\partial t}
% \,,\quad Q(s) = u_{01}(s,s,x_0)
% \\
%\end{eqnarray*}
%
 An insulated equation for 
$Q(t)$ 
 is this: 

\begin{eqnarray*}
\makebox[18ex][l]{$\displaystyle 
 \Big(\frac{\partial^2 }{\partial t^2}
 +{2\gamma}\frac{\partial }{\partial t}
 + \Omega^2 \Big) 
 \Big(Q(t)-u_{01}(t,s,x_0)\Big)
$}%%%end of makebox 
\\
 &=&
 {2\gamma}\frac{\partial q(t)}{\partial t}\Big|_{t=s} 
 - 
 {2\gamma}\frac{\partial  u_{01}(t,s,x_0)}{\partial t}
 + {2\gamma} \int_s^t f_0(\tau) d\tau
 \\
\end{eqnarray*}
 The proper initial relations are these: 
$$
  \Big(Q(s)-u_{01}(s,s,x_0)\Big) = 0\,,\,
  \frac{\partial\Big(Q(t)-u_{01}(t,s,x_0)\Big)}{\partial t}\Big|_{t=s}
 = {2\gamma}\Big(q(s)-u_{01}(s,s,x_0)\Big) 
$$%
 For  
$u=u(t,x)$,
 we have concluded:
\begin{eqnarray*}
 u(t,x) 
 &=&
 -{2\gamma}
 \int_s^t \Big( 1_{+}(x+c(t-\tau)-x_0)-1_{+}(x-c(t-\tau)-x_0)\Big)
 \Big(Q(\tau)-q(\tau)\Big)d\tau
 \\&&{} + u_{01}(t,s,x) 
\end{eqnarray*}

 Other expressions for 
$u=u(t,x)$ 
 are these:
\begin{eqnarray*}
 u(t,x) 
 &=&
 \left\{ 
 \begin{array}{ccl}
 \displaystyle
 Q(t-|x-x_0|/c)-u_{01}(t-|x-x_0|/c,s,x_0)  
                                  &,& \mbox{ if } s \leq t-|x-x_0|/c \\ 
                 0                &,& \mbox{ if } t-|x-x_0|/c < s \leq t \\ 
                 0                &,& \mbox{ if } t \leq s < t+|x-x_0|/c \\ 
 -Q(t+|x-x_0|/c)+u_{01}(t+|x-x_0|/c,s,x_0) 
                                  &,& \mbox{ if } t+|x-x_0|/c \leq s \\ 
 \end{array}
 \right\} 
 \\[\medskipamount]\\
 &&{}\qquad + u_{01}(t,s,x) 
\end{eqnarray*}
 where 
$u_{01}$ 
 is defined by the relations 
\begin{eqnarray*}
 u_{0}(t,s,x)
 &:=& 
 c_{+}(x+ct) + c_{-}(x-ct)
 \\&=&
 \frac12 \Big(u(s,x+c(t-s)) + u(s,x-c(t-s))\Big) 
   + \frac{1}{2c}\int_{x-c(t-s)}^{x+c(t-s)}
 (\frac{\partial u(t,\xi)}{\partial t})\Big|_{t=s} d\xi 
 \\
%\end{eqnarray*}
%
%
%\begin{eqnarray*}
 u_{01}(t,s,x)
 &:=& 
  \frac{1}{2c}\int_s^t 
  \Big(\tilde f_1(\tau,x+c(t-\tau))-\tilde f_1(\tau,x-c(t-\tau))\Big) d\tau 
  +  u_{0}(t,s,x)
%c_{+}(x+ct) + c_{-}(x-ct)
\end{eqnarray*}

\newpage 
 In the following subsections, we will mostly discuss the case where 
$$
 c=1\,, s=0 \,, x_0 = 0 \,, f_1 =0 \,. 
$$
 it will be convenient to have rewritten some of the recent formulae 
 in the proper way.

$$
 \Bigl(\frac{\partial^2 }{\partial t^2} 
  +{2\gamma}\frac{\partial }{\partial t} 
  +\Omega^2 \Bigr)q(t) 
  =
   {2\gamma}\frac{\partial q(t)}{\partial t}\Big|_{t=0} 
   + \Omega^2 u_{0}(t,0,0) 
   + f_0(t) 
   +{2\gamma} \int_0^t f_0(\tau) d\tau 
$$ 

\begin{eqnarray*}
 Q(t)
 &=&
 -{2\gamma}
 \int_0^t 
 \Big(Q(\tau)-q(\tau)\Big) d\tau 
 + u_{0}(t,0,0)
 \\
\end{eqnarray*}
 An insulated equation for 
$Q(t)$ 
 is this: 
\begin{eqnarray*} 
 \Big(\frac{\partial^2 }{\partial t^2}
 +{2\gamma}\frac{\partial }{\partial t}
 + \Omega^2 \Big) 
 \Big(Q(t)-u_{0}(t,0,0)\Big)
 &=&
 {2\gamma}\frac{\partial q(t)}{\partial t}\Big|_{t=0} 
 - 
 {2\gamma}\frac{\partial  u_{0}(t,0,0)}{\partial t}
 + {2\gamma} \int_0^t f_0(\tau) d\tau
 \\
\end{eqnarray*}
 The proper initial conditions are these: 
$$
  \Big(Q(0)-u_{0}(0,0,0)\Big) = 0\,,\,
  \frac{\partial\Big(Q(t)-u_{0}(t,0,0)\Big)}{\partial t}\Big|_{t=0}
 = {2\gamma}\Big(q(0)-u_{0}(0,0,0)\Big) 
$$%
 Expressions for  
$u=u(t,x)$,
 are these:
\begin{eqnarray*}
 u(t,x) 
 &=&
 \left\{ 
 \begin{array}{ccl}
 \displaystyle
 Q(t-|x|)-u_{0}(t-|x|,0,0)  &,& \mbox{ if } 0 \leq t-|x| \\ 
                 0                &,& \mbox{ if } t-|x| < 0 \leq t \\ 
                 0                &,& \mbox{ if } t \leq 0 < t+|x| \\ 
 -Q(t+|x|)+u_{0}(t+|x|,0,0) &,& \mbox{ if } t+|x| \leq 0 \\ 
 \end{array}
 \right\} 
 \\[\medskipamount]\\
 &&{}\qquad + u_{0}(t,0,x) 
\end{eqnarray*}
 where 

\begin{eqnarray*}
 u_{0}(t,0,x)
 &:=& 
 c_{+}(x+t) + c_{-}(x-t)
 \\&=&
 \frac12 \Big(u(0,x+t) + u(0,x-t)\Big) 
   + \frac{1}{2}\int_{x-t}^{x+t}
 (\frac{\partial u(t,\xi)}{\partial t})\Big|_{t=0} d\xi 
\end{eqnarray*}

\newpage 
\subsection
{ A Standard Model. D'Alembert-Kirchhoff-like formulae }

 The considered in the previous subsection model is not standard, 
 if one means classical or quantum field models.
 In that field, a standard looks rather like this:
\begin{eqnarray*}
 \frac{\partial^2 q_0(t)}{\partial t^2}
 &=&-\Omega^2 q_0(t) + \gamma_1 Q_{\phi}(t) + f_0(t) 
\\
 \frac{\partial^2 \phi(t,x)}{\partial t^2} 
 &=&
 c^2\frac{\partial^2 \phi(t,x)}{\partial x^2}
 +4\gamma_2 c\delta(x-x_0)q_0(t)+f_1(t,x) 
\\
 Q_{\phi}(t)  
 &=& \phi(t,x_0)   
\\
\end{eqnarray*}
 In the latter situation one can obtain that 
\begin{eqnarray*}
 Q_{\phi}(t)
 &=&
 {2\gamma_2}
 \int_s^t 
 q_0(\tau) d\tau 
 + \phi_{01}(t,s,x_0)
 \\
\end{eqnarray*}
 for the proper 
%%  suited, adapted, becoming, correct 
$\phi_{01}(t,x)$\,.
  Next, with the reasons of the previous subsections: 
$$
\fbox{
\parbox{0.95\textwidth}{
\begin{eqnarray*}
 \frac{\partial^2 q_0(t)}{\partial t^2}
 &=&-\Omega^2 q_0(t) + 
 {2\gamma_2}\gamma_1 
 \int_s^t 
 q_0(\tau) d\tau 
 + \gamma_1 \phi_{01}(t,s,x_0)
 + f_0(t) 
\end{eqnarray*}
 }%%end of parbox
    }%%end of fbox
$$
 and   
\begin{eqnarray*}
\makebox[18ex][l]{$\displaystyle 
 \Big(\frac{\partial^2 }{\partial t^2} + \Omega^2 \Big)
 \Big(Q_{\phi}(t)- \phi_{01}(t,s,x_0)\Big)
$}%%%end of makebox 
\\
 &=&
 {2\gamma_2}
 \Big(\frac{\partial^2 }{\partial t^2} + \Omega^2 \Big)
 \int_s^t 
 q_0(\tau) d\tau 
 \\
 &=&
 {2\gamma_2}
 \int_s^t 
 \Big(\frac{\partial^2 }{\partial \tau^2} + \Omega^2 \Big)
 q_0(\tau) d\tau 
  +{2\gamma_2}\frac{\partial }{\partial t}q(t)\Big|_{t=s}
 \\
 &=&
 {2\gamma_2}
 \int_s^t 
 \Big(\gamma_1 Q_{\phi}(\tau) + f_0(\tau)\Big) d\tau
  +{2\gamma_2}\frac{\partial }{\partial t}q(t)\Big|_{t=s}
 \\
 &=&
 {2\gamma_2}\gamma_1
 \int_s^t 
 \Big( Q_{\phi}(\tau) - \phi_{01}(\tau,s,x_0)\Big) d\tau
\\&&{ }
 +  
 {2\gamma_2}
 \int_s^t 
 \Big( \gamma_1\phi_{01}(\tau,s,x_0)+ f_0(\tau)\Big) d\tau
  +{2\gamma_2}\frac{\partial }{\partial t}q(t)\Big|_{t=s}
\end{eqnarray*}
 Thus 
$$
\fbox{
\parbox{0.95\textwidth}{
\begin{eqnarray*}
\makebox[18ex][l]{$\displaystyle 
 \Big(\frac{\partial^2 }{\partial t^2} + \Omega^2 \Big)
 \Big(Q_{\phi}(t)- \phi_{01}(t,s,x_0)\Big)
$}%%%end of makebox 
\\
 &=&
 {2\gamma_2}\gamma_1
 \int_s^t 
 \Big( Q_{\phi}(\tau) - \phi_{01}(\tau,s,x_0)\Big) d\tau
\\&&{ }
 +  
 {2\gamma_2}
 \int_s^t 
 \Big( \gamma_1\phi_{01}(\tau,s,x_0)+ f_0(\tau)\Big) d\tau
  +{2\gamma_2}\frac{\partial }{\partial t}q(t)\Big|_{t=s}
\end{eqnarray*}
 }%%end of parbox
    }%%end of fbox
$$

 We can reduce these equations to {\bf ordinary } differential equations, but  
 these ones become {\bf third } order equations,   
\begin{eqnarray*}
 \Big(\frac{\partial^3 }{\partial t^3} 
  + \Omega^2\frac{\partial }{\partial t}
  - {2\gamma_2}\gamma_1 
 \Big)
 q_0(t)
 &=&
 \gamma_1
 \frac{\partial \phi_{01}(t,s,x_0)}{\partial t}
 + \frac{\partial f_0(t)}{\partial t}
 \\
 \Big(\frac{\partial^3 }{\partial t^3} 
  + \Omega^2\frac{\partial }{\partial t}
  + {2\gamma_2}\gamma_1 
 \Big)
 \Big(Q_{\phi}(t)- \phi_{01}(t,s,x_0)\Big)
 &=&
 {2\gamma_2}
 \Big( \gamma_1\phi_{01}(t,s,x_0)+ f_0(t)\Big) 
\end{eqnarray*}
 subject to the supplementary  initial relations
\begin{eqnarray*}
 \Big(\frac{\partial^2 }{\partial t^2} + \Omega^2 \Big)
 q_0(t) 
 \Big|_{t=s} 
  &=& f_0(s) 
 \\
 \Big(\frac{\partial^2 }{\partial t^2} + \Omega^2 \Big)
 \Big(Q_{\phi}(t)- \phi_{01}(t,s,x_0)\Big)
 \Big|_{t=s} 
  &=&
  {2\gamma_2}\frac{\partial }{\partial t}q(t)\Big|_{t=s}
\end{eqnarray*}
 As for 
$\phi(t,x)$, 
 one can obtain that 
\begin{eqnarray*}
 \phi(t,x) 
 &=&
 {2\gamma_2}
 \int_s^t \Big( 1_{+}(x+c(t-\tau)-x_0)-1_{+}(x-c(t-\tau)-x_0)\Big)
 q_0(\tau)d\tau
 \\&&{} + \phi_{01}(t,s,x) 
\end{eqnarray*}
 Other expressions for 
$\phi(t,x)$ 
 will be these:
\begin{eqnarray*}
 \phi(t,x) 
 &=&
 \left\{ 
 \begin{array}{ccl}
 \displaystyle
 {2\gamma_2}
 \int_s^{t-|x-x_0|/c} 
 q_0(\tau)d\tau &,& \mbox{ if } s \leq t-|x-x_0|/c \\ 
       0        &,& \mbox{ if } t-|x-x_0|/c < s \leq t \\ 
       0        &,& \mbox{ if } t \leq s < t+|x-x_0|/c \\ 
 \displaystyle
 -{2\gamma_2}
 \int_s^{t+|x-x_0|/c} 
 q_0(\tau)d\tau &,& \mbox{ if } t+|x-x_0|/c \leq s \\ 
 \end{array}
 \right\} 
 \\[\medskipamount]\\
 &&{}\qquad + \phi_{01}(t,s,x) 
 \\[\bigskipamount] \\ 
 \phi(t,x) 
 &=&
 \left\{ 
 \begin{array}{ccl}
 \displaystyle
 Q_{\phi}(t-|x-x_0|/c)-\phi_{01}(t-|x-x_0|/c,s,x_0)  
                               &,& \mbox{ if } s \leq t-|x-x_0|/c \\ 
                 0             &,& \mbox{ if } t-|x-x_0|/c < s \leq t \\ 
                 0             &,& \mbox{ if } t \leq s < t-|x-x_0|/c \\ 
 -Q_{\phi}(t+|x-x_0|/c)+\phi_{01}(t+|x-x_0|/c,s,x_0) 
                               &,& \mbox{ if } t+|x-x_0|/c \leq s \\ 
 \end{array}
 \right\} 
 \\[\medskipamount]\\
 &&{}\qquad + \phi_{01}(t,s,x) 
\end{eqnarray*}
 where  
$\phi_{01}$
% and  
%$\phi_{0}$\,, 
 is defined 
 of course in the same manner as in the previous subsection,
 by the relations: 
\\
\parbox{\textwidth}{
\begin{eqnarray*}
 \phi_{0}(t,s,x)
 &:=& 
 c_{\phi,+}(x+ct) + c_{\phi,-}(x-ct)
 \\&=&
 \frac12 \Big(\phi(s,x+c(t-s)) + \phi(s,x-c(t-s))\Big) 
   + \frac{1}{2c}\int_{x-c(t-s)}^{x+c(t-s)}
 (\frac{\partial \phi(t,\xi)}{\partial t})\Big|_{t=s} d\xi 
 \\
%\end{eqnarray*}
%
%
%\begin{eqnarray*}
 \phi_{01}(t,s,x)
 &:=& 
   \frac{1}{2c}\int_s^t 
  \Big(\tilde f_1(\tau,x+c(t-\tau))-\tilde f_1(\tau,x-c(t-\tau))\Big) d\tau 
   + \phi_{0}(t,s,x)
%  + c_{\phi,+}(x+ct) + c_{\phi,-}(x-ct)
\end{eqnarray*}
   }%%end of parbox 

 \newpage 

 The case where 
$$
 c=1\,, s=0 \,, x_0 = 0 \,, f_1 =0 \,. 
$$
 is this:

$$
\fbox{
\parbox{0.95\textwidth}{
\begin{eqnarray*}
 \frac{\partial^2 q_0(t)}{\partial t^2}
 &=&-\Omega^2 q_0(t) + 
 {2\gamma_2}\gamma_1 
 \int_0^t 
 q_0(\tau) d\tau 
 + \gamma_1 \phi_{0}(t,0,0)
 + f_0(t) 
\end{eqnarray*}
 }%%end of parbox
    }%%end of fbox
$$
$$
\fbox{
\parbox{0.95\textwidth}{
\begin{eqnarray*}
\makebox[18ex][l]{$\displaystyle 
 \Big(\frac{\partial^2 }{\partial t^2} + \Omega^2 \Big)
 \Big(Q_{\phi}(t)- \phi_{0}(t,0,0)\Big)
$}%%%end of makebox 
\\
 &=&
 {2\gamma_2}\gamma_1
 \int_0^t 
 \Big( Q_{\phi}(\tau) - \phi_{0}(\tau,0,0)\Big) d\tau
\\&&{ }
 +  
 {2\gamma_2}
 \int_0^t 
 \Big( \gamma_1\phi_{0}(\tau,0,0)+ f_0(\tau)\Big) d\tau
  +{2\gamma_2}\frac{\partial }{\partial t}q(t)\Big|_{t=0}
\end{eqnarray*}
 }%%end of parbox
    }%%end of fbox
$$

 Reducing to ordinary differential equations:
\begin{eqnarray*}
 \Big(\frac{\partial^3 }{\partial t^3} 
  + \Omega^2\frac{\partial }{\partial t}
  - {2\gamma_2}\gamma_1 
 \Big)
 q_0(t)
 &=&
 \gamma_1
 \frac{\partial \phi_{0}(t,0,0)}{\partial t}
 + \frac{\partial f_0(t)}{\partial t}
 \\
 \Big(\frac{\partial^3 }{\partial t^3} 
  + \Omega^2\frac{\partial }{\partial t}
  + {2\gamma_2}\gamma_1 
 \Big)
 \Big(Q_{\phi}(t)- \phi_{0}(t,0,0)\Big)
 &=&
 {2\gamma_2}
 \Big( \gamma_1\phi_{0}(t,0,0)+ f_0(t)\Big) 
\end{eqnarray*}
 supplementary initial relations:
\begin{eqnarray*}
 \Big(\frac{\partial^2 }{\partial t^2} + \Omega^2 \Big)
 q_0(t) 
 \Big|_{t=0} 
  &=& f_0(0) 
 \\
 \Big(\frac{\partial^2 }{\partial t^2} + \Omega^2 \Big)
 \Big(Q_{\phi}(t)- \phi_{0}(t,0,0)\Big)
 \Big|_{t=0} 
  &=&
  {2\gamma_2}\frac{\partial }{\partial t}q(t)\Big|_{t=0}
\end{eqnarray*}
 As for  
$\phi(t,x)$, 

\begin{eqnarray*}
 \phi(t,x) 
 &=&
 \left\{ 
 \begin{array}{ccl}
 \displaystyle
 Q_{\phi}(t-|x|)-\phi_{0}(t-|x|,0,0)  
                               &,& \mbox{ if } 0 \leq t-|x| \\ 
                 0             &,& \mbox{ if } t-|x| < 0 \leq t \\ 
                 0             &,& \mbox{ if } t \leq 0 < t+|x| \\ 
 -Q_{\phi}(t+|x|)+\phi_{0}(t+|x|,0,0) 
                               &,& \mbox{ if } t+|x| \leq 0 \\ 
 \end{array}
 \right\} 
 \\[\medskipamount]\\
 &&{}\qquad + \phi_{0}(t,0,x) 
\end{eqnarray*}
 where  
$\phi_{0}$
% and  
%$\phi_{0}$\,, 
 is defined 
 by the relations: 
\\
\parbox{\textwidth}{
\begin{eqnarray*}
 \phi_{0}(t,0,x)
 &:=& 
 c_{\phi,+}(x+t) + c_{\phi,-}(x-t)
 \\&=&
 \frac12 \Big(\phi(0,x+t) + \phi(0,x-t)\Big) 
   + \frac{1}{2}\int_{x-t}^{x+t}
 (\frac{\partial \phi(t,\xi)}{\partial t})\Big|_{t=0} d\xi 
 \\
\end{eqnarray*}
   }%%end of parbox 

  \newpage

\newpage\subsection%*
{ Particular Cases. 1. Radiation Reaction, Braking Radiation } 

 We are all familiar with the fact that 
 any solution to any linear inhomogeneous equation, whatever its nature,  
 is a sum of a solution 
 to the associated homogeneous equation plus 
 arbitrarily taken and fixed solution to 
 the former linear inhomogeneous equation,
 isn't it? 
 Just let us now consider these two cases separately, 
 the case where equations are homogeneous and then inhomogeneous. 

 We start out emphasising that the homogeneous equation array, connected
 to  
$$
 \Bigl(\frac{\partial^2 }{\partial t^2} 
  +{2\gamma}\frac{\partial }{\partial t} 
  +\Omega^2 \Bigr)q(t) 
  =
   {2\gamma}\frac{\partial q(t)}{\partial t}\Big|_{t=0} 
   + \Omega^2 u_{0}(t,0,0) 
   + f_0(t) 
   +{2\gamma} \int_0^t f_0(\tau) d\tau 
$$ 
\begin{eqnarray*} 
 \Big(\frac{\partial^2 }{\partial t^2}
 +{2\gamma}\frac{\partial }{\partial t}
 + \Omega^2 \Big) 
 \Big(Q(t)-u_{0}(t,0,0)\Big)
 &=&
 {2\gamma}\frac{\partial q(t)}{\partial t}\Big|_{t=0} 
 - 
 {2\gamma}\frac{\partial  u_{0}(t,0,0)}{\partial t}
 + {2\gamma} \int_0^t f_0(\tau) d\tau
 \\
\end{eqnarray*}
 is exactly  
\begin{eqnarray*} 
 \Bigl(\frac{\partial^2 }{\partial t^2} 
  +{2\gamma}\frac{\partial }{\partial t} 
  +\Omega^2 \Bigr)q(t) 
  &=&
   {2\gamma}\frac{\partial q(t)}{\partial t}\Big|_{t=0} 
 \\
 \Big(\frac{\partial^2 }{\partial t^2}
 +{2\gamma}\frac{\partial }{\partial t}
 + \Omega^2 \Big) 
 \Big(Q(t)-u_{0}(t,0,0)\Big)
 &=&
 {2\gamma}\frac{\partial q(t)}{\partial t}\Big|_{t=0} 
\end{eqnarray*}
 and NOT 
\begin{eqnarray*} 
 \Bigl(\frac{\partial^2 }{\partial t^2} 
  +{2\gamma}\frac{\partial }{\partial t} 
  +\Omega^2 \Bigr)q(t) 
 &=& 
 0 
 \\ 
 \Big(\frac{\partial^2 }{\partial t^2}
 +{2\gamma}\frac{\partial }{\partial t}
 + \Omega^2 \Big) 
 \Big(Q(t)-u_{0}(t,0,0)\Big)
 &=& 
  0 
\end{eqnarray*}
 The difference is a {\bf rank one } term 
${2\gamma}\frac{\partial q(t)}{\partial t}\Big|_{t=0}$ 
\footnote{ 
 This term, as a function of 
$t$ ,
 is a fixed function of 
$t$, 
e.g. 
$1$, multiplied by a CONSTANT depended on $q$, i.e.,
 by a fixed functional of 
$q$. 
 Using the Dirac's syntax, this term 
 can be written as 
$|a><b|$ 
 with 
$|a>=1$
 and
$<b|q>={2\gamma}\frac{\partial q(t)}{\partial t}\Big|_{t=0}$ . 
 }% end footnote
 . 
 This detail allows us to apply usual machinery of finite rank perturbations 
 theory. 
 Thus, having put 
$$
 \Omega_{\gamma} := \sqrt{\Omega^2 -{\gamma}^2}
$$
 and having taken into account the supplementary initial relations: 
$$
  Q(0) = 0\,,\,
  \frac{\partial Q(t)}{\partial t}\Big|_{t=0}
 = {2\gamma} q(0) 
$$%
 one can show that  
\begin{eqnarray*}
 q(t) 
 &=& 
   e^{-\gamma t} 
  \Big(cos(\Omega_{\gamma} t) 
        +\gamma \frac{sin(\Omega_{\gamma} t)}{\Omega_{\gamma}}
  \Big)
  \Big(q(0) 
      -\frac{2\gamma}{\Omega^2}\frac{\partial q(t)}{\partial t}\Big|_{t=0}  
  \Big)
%%% \\&&{ }
   + e^{-\gamma t}\frac{sin(\Omega_{\gamma} t)}{\Omega_{\gamma}} 
            \frac{\partial q(t)}{\partial t}\Big|_{t=0}  
 \\&&{ }
   +\frac{2\gamma}{\Omega^2}\frac{\partial q(t)}{\partial t}\Big|_{t=0}  
 \\[2\bigskipamount]
%\end{eqnarray*}
%\begin{eqnarray*}
 Q(t) 
 &=& 
   e^{-\gamma t} 
  \Big(cos(\Omega_{\gamma} t) 
        +\gamma \frac{sin(\Omega_{\gamma} t)}{\Omega_{\gamma}}
  \Big)
  \Big( 
      -\frac{2\gamma}{\Omega^2}\frac{\partial q(t)}{\partial t}\Big|_{t=0}  
  \Big)
%%% \\&&{ }
   + e^{-\gamma t}\frac{sin(\Omega_{\gamma} t)}{\Omega_{\gamma}} 
          2 \gamma q(0)  
 \\&&{ }
   +\frac{2\gamma}{\Omega^2}\frac{\partial q(t)}{\partial t}\Big|_{t=0}  
 \\[2\bigskipamount]
\end{eqnarray*}
 and then 
\begin{eqnarray*}
 u(t,x)  
 &=& 
   e^{-\gamma (t-|x|)}
  \Big(cos\Omega_{\gamma} (t-|x|) 
        +\gamma \frac{sin\Omega_{\gamma}(t-|x|)}{\Omega_{\gamma}}
  \Big)
  \Big( 
      -\frac{2\gamma}{\Omega^2}\frac{\partial q(t)}{\partial t}\Big|_{t=0}  
  \Big)
 \\&&{ }
   + e^{-\gamma (t-|x|)}
   \frac{sin\Omega_{\gamma}(t-|x|)}{\Omega_{\gamma}} 
          2 \gamma q(0)  
 \\&&{ }
   +\frac{2\gamma}{\Omega^2}\frac{\partial q(t)}{\partial t}\Big|_{t=0}  
 \,, \mbox{ if } 0 \leq t-|x| 
 \\
 u(t,x)  
 &=& 
 0
 \,, \mbox{ if } 0 > t-|x| 
 \\
\end{eqnarray*} 
  Notice, 
\begin{eqnarray*} 
 q(t) - Q(t)
 &=&
   e^{-\gamma t}
     cos(\Omega_{\gamma} t)
       q(0)
   + 
       e^{-\gamma t}
         \frac{sin\Omega_{\gamma} t}{\Omega_{\gamma}} 
            \Bigl(
                  \frac{\partial q(t)}{\partial t}\Big|_{t=0} 
                    -2 \gamma q(0)  
            \Bigr)
% - 
% + {2\gamma} \int_0^t f_0(\tau) d\tau
 \\
\end{eqnarray*}
 and then 
$$
 q(t) - Q(t)\to 0 
 \mbox{ as  } t\to +\infty\,,\, \gamma >0
$$ 

  If we concentrate now on the system's behaviour at large  
$t$ , 
 a mathematical detail calls attention. 
 We observe:
 $$
 q(t)\to  {2\gamma}\frac{\partial q(t)}{\partial t}\Big|_{t=0}   
 \mbox{ as } t \to +\infty \,. 
$$  
 Nevertheless, the limit function 
$$
 q_{\infty}(t)={2\gamma}\frac{\partial q(t)}{\partial t}\Big|_{t=0} 
$$
 is NO solution to 
$$
 \frac{\partial^2 q(t)}{\partial t^2} 
  =-{2\gamma}\frac{\partial q(t)}{\partial t} 
   -\Omega^2 q(t) 
   +{2\gamma}\frac{\partial q(t)}{\partial t}\Big|_{t=0} 
$$ 
 every time that 
$$
  {2\gamma}\frac{\partial q(t)}{\partial t}\Big|_{t=0} \not= 0\,, 
$$
 because 
 the 
$q=q_{\infty}$ 
 is a constant, hence its derivative is zero 
 but we have need of  
$$
 \frac{\partial q(t)}{\partial t}\Big|_{t=0} 
 = \frac{\partial q_{\infty}}{\partial t} 
 = \frac{\partial }{\partial t} 
   \Bigl( 
   {2\gamma}\frac{\partial q(t)}{\partial t}\Big|_{t=0} 
   \Bigr) = 0 \,. 
$$
 Similar phenomena, 
 one can detect them in the electrodynamics of moving charges.
 
\addvspace{\bigskipamount}

 Now, we are directing our attention to the fact that the oscillator 
 moves with damped amplitude:
$$
 |q(t)-q_{\infty}| \leq e^{-\gamma t}\cdot const \,.
$$ 

 The reason is very plain. 
 When the particle (oscillator) begins to move, it makes 
 some region of the field `move': 
 the oscillator generates, emits waves, 
 --one is used to saying: `the oscillator radiates'. 

 Now then, the oscillator radiates.
 It entails some energy expenses
 and thus the field, made move, 
 eventually makes the conduct of the oscillator change.
 This phenomenon is said to be 
 a {\bf radiation reaction}. 

\addvspace{\bigskipamount}
 
 Next, since the oscillator has emitted waves, 
 waves run away and brings away some 
 portions of the oscillator energy,
 braking the moving of the oscillator. 
 If something had reflected the waves, more precise, 
 if the emitted waves had returned to the oscillator, 
 we could wait for that the these waves would stimulate an increase 
 of the amplitude of the oscillator. 
 But if no wave returns to the oscilator, 
 --it is exactly our case--, 
 the damping dominates and the amplitude decreases 
\footnote{
 as we have seen, it decreases to the zero value
}% end footnote 
 . 
 This phenomenon is said to be a
 {\bf braking radiation} or {\bf damping radiation}.

\newpage\subsection%*
{ Particular Cases. 2. Resonance } 

 We begin to analyse the inhomogeneous equations with the most simple 
 case, where there is an only right-hand incident wave   
$$
  u_{01}(t,0,x) =  u_{0}(t,0,x) = c_{+}(x+t) 
$$
 and where one can easy compute the solution: 
 don't deviate from standards, let
$$
  u_{01}(t,0,x) =  u_{0}(t,0,x) = c_{+}(x+t) 
  := A\sin(k(x+t))
$$
 i.e.,  
\begin{eqnarray*} 
 \Bigl(\frac{\partial^2 }{\partial t^2} 
  +{2\gamma}\frac{\partial }{\partial t} 
  +\Omega^2 \Bigr)q(t) 
  &=&
   {2\gamma}\frac{\partial q(t)}{\partial t}\Big|_{t=0} 
   + \Omega^2 A\sin(kt)
\end{eqnarray*}
\begin{eqnarray*} 
 \Big(\frac{\partial^2 }{\partial t^2}
 +{2\gamma}\frac{\partial }{\partial t}
 + \Omega^2 \Big) 
 \Big(Q(t)-A\sin(kt)\Big)
 &=&
 {2\gamma}\frac{\partial q(t)}{\partial t}\Big|_{t=0} 
 - 
 {2\gamma}Ak\cos(kt)
 \\
\end{eqnarray*}
 According to the standard reading of the universiy course to 
 the theory of ordinary differential equations, if 
\begin{eqnarray*} 
 \Bigl(\frac{\partial^2 }{\partial t^2} 
  +{2\gamma}\frac{\partial }{\partial t} 
  +\Omega^2 \Bigr)y(t) 
  &=&
   const 
   + \Omega^2 A\sin(kt)
\end{eqnarray*}
\begin{eqnarray*} 
 \Big(\frac{\partial^2 }{\partial t^2}
 +{2\gamma}\frac{\partial }{\partial t}
 + \Omega^2 \Big) 
 Y(t)
 &=&
   const 
 - 
 {2\gamma}Ak\cos(kt)
 \\
\end{eqnarray*}
 then a partial solution to these equations is given by: 
\begin{eqnarray*}
 y(t)
 &=& 
 \frac{\Omega^2 A}{ (-k^2+\Omega^2)^2 + (2\gamma k)^2 } 
  \left( 
    ( -k^2+\Omega^2 )\sin(kt) -2\gamma k \cos(kt)
  \right) 
 +
 \frac{const}{\Omega^2}
\end{eqnarray*}
\begin{eqnarray*}
 Y(t)
 &=& 
 \frac{-2\gamma Ak}{ (-k^2+\Omega^2)^2 + (2\gamma k)^2 } 
  \left( 
     2\gamma k\sin(kt) + ( -k^2+\Omega^2 ) \cos(kt)
  \right) 
 +
 \frac{const}{\Omega^2}
\end{eqnarray*}

 On the surface, there is no essential deflecting from the situation in the 
 university course of ordinary differential equation, but it is 
 only on the surface.   

 Through the presence of the term 
$\frac{\partial q(t)}{\partial t}\Big|_{t=0}$ 
 the forms 
\begin{eqnarray*}
 q(t)
 &=& 
    A_s\sin(kt) + A_c\cos(kt) + Const 
\end{eqnarray*}
 need not contain a particular solution to the former equation.
 Although, we have, of course,   
\begin{eqnarray*}
 q(t)
 &=& 
 \frac{\Omega^2 A}{ (-k^2+\Omega^2)^2 + (2\gamma k)^2 } 
  \left( 
    ( -k^2+\Omega^2 )\sin(kt) -2\gamma k \cos(kt)
  \right) + const_1
\\&&{} 
 +
 \frac{2\gamma}{\Omega^2}\frac{\partial q(t)}{\partial t}\Big|_{t=0}
\\&&{} 
 +
   e^{-\gamma t} 
  \Bigl(const_1\cos(\Omega_{\gamma} t) 
        +const_2 \sin(\Omega_{\gamma} t)
  \Bigr) \,;
\end{eqnarray*}
 we can write this relation as follows: 
\begin{eqnarray*}
 q(t)
 &=& 
 \frac{\Omega^2 A}{\sqrt{ (-k^2+\Omega^2)^2 + (2\gamma k)^2}} 
   \sin(kt +\phi_k )
 + const_1
\\&&{} 
 +
 \frac{2\gamma}{\Omega^2}\frac{\partial q(t)}{\partial t}\Big|_{t=0}
\\&&{} 
 +
   e^{-\gamma t} 
  \Bigl(const_1\cos(\Omega_{\gamma} t) 
        +const_2 \sin(\Omega_{\gamma} t)
  \Bigr)
\end{eqnarray*}
 where
$\phi_k$ 
 is such that 
$$
  \cos\phi_k = \frac{-k^2+\Omega^2}{\sqrt{(-k^2+\Omega^2)^2 + (2\gamma k)^2 }} 
\,,\,
  \sin\phi_k = -\frac{2\gamma k}{\sqrt{(-k^2+\Omega^2)^2 + (2\gamma k)^2 }} 
\,. 
$$

\addvspace{\bigskipamount}

 As for 
$Q(t)-A\sin(kt)$, 
 on referring to the equation to 
$Q(t)-A\sin(kt)$ 
 we see, on the contrary, that the term 
${2\gamma}\frac{\partial q(t)}{\partial t}\Big|_{t=0}$ 
 is there exterior and does not hinder for immediate displaying 
 a particular solution to the {lonely} equation to 
$Q(t)-A\sin(kt)$ 
 itself. 
 It is this: 
\begin{eqnarray*}
 Q(t)-A\sin(kt)
 &=& 
 \frac{-2\gamma Ak}{ (-k^2+\Omega^2)^2 + (2\gamma k)^2 } 
  \left( 
     2\gamma k\sin(kt) + ( -k^2+\Omega^2 ) \cos(kt)
  \right) 
 +
 \frac{2\gamma}{\Omega^2}\frac{\partial q(t)}{\partial t}\Big|_{t=0}
\end{eqnarray*}
 Nevertheless, this form is improper, because of the initial relations 
$$
  \Big(Q(0)-c_{+}(0)\Big) = 0\,,\,
  \frac{\partial\Big(Q(t)-c_{+}(t)\Big)}{\partial t}\Big|_{t=0}
 = {2\gamma}\Big(q(0)-c_{+}(0)\Big) 
$$%
 which in our case read 
$$
  Q(0)= 0\,,\,
  \frac{\partial\Big(Q(t)-A\sin(kt)\Big)}{\partial t}\Big|_{t=0}
 = {2\gamma}q(0) 
$$
 Once again, we need to take into account  
 the terms containing 
$e^{-\gamma t}\Bigl(\cdots\Bigl)$ 
 : 
\begin{eqnarray*}
 Q(t)-A\sin(kt)
 &=& 
 \frac{-2\gamma Ak}{ (-k^2+\Omega^2)^2 + (2\gamma k)^2 } 
  \left( 
     2\gamma k\sin(kt) + ( -k^2+\Omega^2 ) \cos(kt)
  \right) 
 +
 \frac{2\gamma}{\Omega^2}\frac{\partial q(t)}{\partial t}\Big|_{t=0}
 \\&&{}
  +e^{-\gamma t} 
  \Bigl(const_3\cos(\Omega_{\gamma} t) 
        +const_4 \sin(\Omega_{\gamma} t)
  \Bigr)
\\
 &=& 
 \frac{-2\gamma Ak}{\sqrt{ (-k^2+\Omega^2)^2 + (2\gamma k)^2}} 
   \cos(kt +\phi_k )
 +
 \frac{2\gamma}{\Omega^2}\frac{\partial q(t)}{\partial t}\Big|_{t=0}
 \\&&{}
  +e^{-\gamma t} 
  \Bigl(const_3\cos(\Omega_{\gamma} t) 
        +const_4 \sin(\Omega_{\gamma} t)
  \Bigr) 
\end{eqnarray*}

\addvspace{\bigskipamount}

 Let us restrict ourselves to the case where 
$\gamma t >> 1$.  
 In this case 
$e^{-\gamma t}\Bigl(\cdots\Bigl)\approx 0$ 
 and 
\begin{eqnarray*}
 q(t)
 &\approx & 
 \frac{\Omega^2 A}{\sqrt{ (-k^2+\Omega^2)^2 + (2\gamma k)^2}} 
   \sin(kt +\phi_k )
 + const_1
\\&&{} 
 +
 \frac{2\gamma}{\Omega^2}\frac{\partial q(t)}{\partial t}\Big|_{t=0}
\end{eqnarray*}
\begin{eqnarray*}
 Q(t)-A\sin(kt)
 &\approx & 
 \frac{-2\gamma Ak}{ (-k^2+\Omega^2)^2 + (2\gamma k)^2 } 
  \left( 
     2\gamma k\sin(kt) + ( -k^2+\Omega^2 ) \cos(kt)
  \right) 
 +
 \frac{2\gamma}{\Omega^2}\frac{\partial q(t)}{\partial t}\Big|_{t=0} \,. 
\end{eqnarray*}
  If 
$k^2=\Omega^2 - 2\gamma^2$,
\footnote{ fix on: $k^2=\Omega^2 - 2\gamma^2$, not $k^2=\Omega^2$ !! } 
 then the amplitude of the harmonic part of 
$q(t)$ , 
 i.e., the value of the quantity  
$$
 \frac{|\Omega^2 A|}{\sqrt{ (-k^2+\Omega^2)^2 + (2\gamma k)^2}} \,, 
$$
 becomes maximal
\footnote{ up to assumption made just now }: 
 This phenomenon is said to be a {\bf resonance }.

\addvspace{\bigskipamount}

 Next, if 
$k^2=\Omega^2$,
\footnote{ fix on: $k^2=\Omega^2$, not $k^2=\Omega^2 - 2\gamma^2$ !! } 
 and if, in addition, 
$\frac{\partial q(t)}{\partial t}\Big|_{t=0}=0$,
 then 
\begin{eqnarray*}
 Q(t) - A\sin(kt)
 &\approx & 
 \frac{-2\gamma Ak}{(2\gamma k)^2 } 
  \left( 
     2\gamma k\sin(kt)
  \right) 
\\
  &=&{ } 
  -A\sin(kt) 
\end{eqnarray*}
\footnote{ 
 hence, 
$Q(t)\approx 0$}

 Let us now recall, that 
\begin{eqnarray*}
 u(t,x) 
 &=&
 \left\{ 
 \begin{array}{ccl}
 \displaystyle
 Q(t-|x|)-c_{+}(t-|x|)          &,& \mbox{ if } 0 \leq t-|x| \\ 
               0                &,& \mbox{ if } t-|x| < 0 \leq t \\ 
               0                &,& \mbox{ if } t \leq 0 < t+|x| \\ 
 -Q(t+|x|)+c_{+}(t+|x|)         &,& \mbox{ if } t+|x| \leq 0 \\ 
 \end{array}
 \right\} 
% \\[\medskipamount]\\
% &&{}\qquad 
 + c_{+}(x+t) 
\end{eqnarray*}
 and hence, if 
$k^2=\Omega^2$,  
$\frac{\partial q(t)}{\partial t}\Big|_{t=0}=0$,
$\gamma t >> 1$,  
 then 
\begin{eqnarray*}
 u(t,x) 
 &\approx &
 \left\{ 
 \begin{array}{ccl}
 \displaystyle
 -A\sin(k(t-|x|))          &,& \mbox{ if } 0 \leq t-|x| \\ 
               0                &,& \mbox{ if } t-|x| < 0 \leq t \\ 
               0                &,& \mbox{ if } t \leq 0 < t+|x| \\ 
 -A\sin(k(t+|x|))         &,& \mbox{ if } t+|x| \leq 0 \\ 
 \end{array}
 \right\} 
% \\[\medskipamount]\\
% &&{}\qquad 
 + A\sin(k(x+t)) 
\end{eqnarray*}
 In particular,
\begin{eqnarray*}
 u(t,x) 
 &\approx &
 \left\{ 
 \begin{array}{ccl}
 \displaystyle
  0                  &,& \mbox{ if } 0 \leq t-|x| , x< 0,\gamma t >> 1  
\\ 
  
  2A\sin(kx)cos(kt)  &,& \mbox{ if } 0 \leq t-|x| , x> 0, \gamma t >> 1
\\ 
    A\sin(k(x+t))    &,& \mbox{ if } t-|x| < 0 \leq t \\ 
 \end{array}  
 \right\} 
\end{eqnarray*}

\addvspace{\bigskipamount}

 We see, in the case of 
$$
 k^2=\Omega^2 \,,\,
 \frac{\partial q(t)}{\partial t}\Big|_{t=0}=0  \,,\,
 \gamma t >> 1 
$$
 the incident wave is essentially completely reflected by the oscillator!!! 
 This phenomenon can also be referred to as a kind of {\bf resonance}.

\newpage\subsection%*%
{ Comments. Discussion.} 

 In this section,  
 we have considered two models of interaction and obtained 
 for oscillators the corresponding equations:
$$
 \frac{\partial^2 q(t)}{\partial t^2} 
  =-{2\gamma}\frac{\partial q(t)}{\partial t} 
   -\Omega^2 q(t) 
   +{2\gamma}\frac{\partial q(t)}{\partial t}\Big|_{t=0} 
   + \Omega^2 u_{0}(t,0,0) 
   + f_0(t) 
   +{2\gamma} \int_0^t f_0(\tau) d\tau 
$$ 

\begin{eqnarray*}
 \frac{\partial^2 q_0(t)}{\partial t^2}
 &=&-\Omega^2 q_0(t) + 
 {2\gamma_2}\gamma_1 
 \int_0^t 
 q_0(\tau) d\tau 
 + \gamma_1 \phi_{01}(t,0,0)
 + f_0(t) 
\end{eqnarray*}
 These equations 
 are not literally {\bf ordinary } differential equations ,
% without deviating from the definitions, 
% literally, verbatim, in the proper sense of the words
% undeviating, undeviatingly   
{\bf ordinary } differential equations,
 because of `unordinary' terms  
$$
   +{2\gamma}\frac{\partial q(t)}{\partial t}\Big|_{t=0}   
 \mbox{ and, resp., } 
 {2\gamma_1\gamma_2}
 \int_0^t 
 q_0(\tau) d\tau 
$$
 We can reduce these equations to {\bf ordinary } differential equations, 
 but the latter turn out to be {\bf third } order equations:
$$
 \frac{\partial^3 q(t)}{\partial t^3}  
  =-{2\gamma}\frac{\partial^2 q(t)}{\partial t^2} 
   -\Omega^2 \frac{\partial q(t)}{\partial t}  
 +\frac{\partial }{\partial t} 
  \Big( \Omega^2 u_{0}(t,0,0)
   + f_0(t) 
   +{2\gamma} \int_0^t f_0(\tau) d\tau \Big) 
$$ 

\begin{eqnarray*}
 \frac{\partial^3 q_0(t)}{\partial t^3}
 &=&-\Omega^2 \frac{\partial q_0(t)}{\partial t} + \gamma_1 
 \Big(
 {2\gamma_2}
 q_0(t) 
 + \frac{\partial u_{0}(t,0,0)}{\partial t}
 \Big) 
 + \frac{\partial f_0(t)}{\partial t} 
\end{eqnarray*}

 We concentrate now on the system's behaviour at large  
$t$ . 
 In order to estimate the asymptotic behaviour of 
$q(t)$ 
 and 
$q_0(t)$,  
 as 
$ t \to +\infty $,
 let us handle with the characteristic polynomials. They are:

$$
 \lambda^3 + {2\gamma}\lambda^2+\Omega^2\lambda  \quad, \quad 
 \lambda_0^3 +\Omega^2\lambda_0 - 2\gamma_1{\gamma_2} 
$$ 
 The first polynomial has the roots
$$
 \lambda_1 = 0,  
 \lambda_2 = -{\gamma} 
             +\sqrt{{\gamma}^2-\Omega^2} , 
 \lambda_3 = -{\gamma} 
             -\sqrt{{\gamma}^2-\Omega^2} , 
$$
 and we see that first root is zero and two other roots
 have non-positive real parts. 
 Hence, except for the unique case, the case of 
$$
   +{2\gamma}\frac{\partial q(t)}{\partial t}\Big|_{t=0} 
   \not= 0 \,,
$$
 we have exponentially decaying 
$q(t)$
 as 
$t \to +\infty$. 
 Whatever the case,  
$q(t)$
 is bounded as 
$t \to +\infty$. 
  Recall, if 
$f_0$ 
 and 
$u_{0}$ 
 both are identically zero, then 
$$
 q(t)\to  {2\gamma}\frac{\partial q(t)}{\partial t}\Big|_{t=0}   
 \mbox{ as } t \to +\infty \,. 
$$  
 However, the limit function 
$$
 q_{\infty}(t)={2\gamma}\frac{\partial q(t)}{\partial t}\Big|_{t=0} 
$$
 is NO solution to 
$$
 \frac{\partial^2 q(t)}{\partial t^2} 
  =-{2\gamma}\frac{\partial q(t)}{\partial t} 
   -\Omega^2 q(t) 
   +{2\gamma}\frac{\partial q(t)}{\partial t}\Big|_{t=0} 
$$ 
 every time that 
$$
  {2\gamma}\frac{\partial q(t)}{\partial t}\Big|_{t=0} \not= 0\,. 
$$

 As for the second polynomial, and the second equation, 
 the situation is more dramatic, much more. The polynomial 
$$
 \lambda_0^3 +\Omega^2\lambda_0 - 2\gamma_1{\gamma_2} 
$$ 
 has one pure real root 
$\lambda_{01}$ 
 and two complex-conjugated ones: 
$\lambda_{02}\,,\lambda_{03} \,,\,\lambda_{03} = \overline{\lambda_{02}}$ .
 Since 
$$
  \lambda_{01} +\lambda_{02} +\lambda_{03}  = 0 
 \,,\,
  \lambda_{01} \lambda_{02} \lambda_{03}  = 2\gamma_1{\gamma_2} 
$$ 
 we have 
$$
  \lambda_{01} +2Re\lambda_{02} = 0 
 \,,\,
  \lambda_{01} |\lambda_{02}|^2  = 2\gamma_1{\gamma_2} 
$$ 
 and hence 
$$
 \lambda_{01} > 0\,\mbox{ (!!!) },
 Re\lambda_{02} = Re\lambda_{03} <0 \,. 
$$
 Thus, we observe an EXPONENTIAL GROWTH of the oscillator amplitude, 
 as 
$t\to +\infty$, 
 the fact used to perplexing physicist's mind. 
 We cannot here hope we have simply confused the `time direction'.
 If we had, 
 we would have {\bf two} roots with strictly positive real parts! 
 We defer the more detailed discussion on this subject
 and notice only, that
 a related phenomenon is known in electrodynimics, see 
 Abraham-Lorentz-Dirac equations.

%%%%%%%%%%%
%%%%%%%%%%%

\newpage\section%*%
{ A little more general Model with Finite Rank Interaction. 
  The Case of non-local Interaction and arbitrary Interaction Forces.} 

 We will now slightly change the model. 
 We will not suppose that 
$\delta $
 shall be Dirac's 
$\delta$-function 
\footnote{ 
 In order to no confusion can ocurr we replace the symbol
$\delta$ 
 by 
$\delta_{\alpha}$ 
 }%% end of this footnote
 and 
$B$ 
 be 
$\frac{\partial^2}{\partial t^2}$
 on the whole line. 
 At the beginning of this section, 
 we will consider an abstract d'Alambert-like equation 
 and recall the associated d'Alembert-like solution formulae.
 Then we introduce an abstract analogue of the interaction 
 discussed in the previous section. 
 Then we discuss some partial cases, imitating the reasons of the 
 same section,  
 and finally indicate how the abstract formulae correlate with 
 the formulae obtained in the previous section.

 Now then. When solving the equation
$$
  \frac{\partial^2 u}{\partial t^2} = Bu+f 
$$
 denote
$$
 u_1:=\frac{\partial u}{\partial t} \,.
$$
 Write the equation to be solved as
$$
 \frac{\partial }{\partial t} 
\Big( \begin{array}{cc}
 u \\ u_1
\end{array}\Big)
 = 
\Big( \begin{array}{cc}
 u_1 \\ Bu + f
\end{array}\Big)
 = 
\Big( \begin{array}{cc}
 0 & 1 \\ B & 0 \\ 
\end{array}\Big)
\Big( \begin{array}{cc}
 u \\ u_1 
\end{array}\Big)
 + 
\Big( \begin{array}{cc}
 0 \\ f 
\end{array}\Big)
$$
 Let
$$
 V_{t,s}
$$
 denote the propagator for 
$$
 \frac{\partial }{\partial t} 
\Big( \begin{array}{cc}
 v \\ v_1 
\end{array}\Big)
 = 
\Big( \begin{array}{cc}
 0 & 1 \\ B & 0 \\
\end{array}\Big)
\Big( \begin{array}{cc}
 v \\ v_1 
\end{array}\Big)
$$
 Then
$$
\Big( \begin{array}{cc}
 u(t) \\ u_1(t) 
\end{array}\Big)
 = V_{t,s}
\Big( \begin{array}{cc}
 u(s) \\ u_1(s) 
\end{array}\Big)
 + \int_s^t V_{t,\tau} 
\Big( \begin{array}{cc}
 0 \\ f(\tau) 
\end{array}\Big)
 d\tau 
$$
 Apply usual matrix -form representation
$$
 V_{t,s} = 
\Big( \begin{array}{cc}
 V_{11}(t,s) & V_{12}(t,s) \\  V_{21}(t,s) & V_{22}(t,s) 
\end{array}\Big)
$$
 Then obtain 
$$
 u(t) = V_{11}(t,s)u(s) 
        + V_{12}(t,s)\Big(\frac{\partial u(t)}{\partial t} \Big|_{t=s} \Big) 
        + \int_s^t V_{12}(t,\tau)f(\tau)d\tau 
$$
$$
 u_{0}(t) := V_{11}(t,s)u(s) 
        + V_{12}(t,s)\Big(\frac{\partial u(t)}{\partial t} \Big|_{t=s} \Big) 
$$
$$
 u(t) = u_{0}(t,s) 
        + \int_s^t V_{12}(t,\tau)f(\tau)d\tau 
$$

 I will now take 
\begin{eqnarray*}
 f(t) 
 &=& -\delta_{\alpha,t,x_0}U_{0}(t,Q) +f_1(t)
\\ 
 Q(t) 
 &=& <l|u(t)> 
\end{eqnarray*}
 where:

 1)  
$U_{0}(t,Q)$ 
 is such that its values must be {\bf numbers}, and 

 2) 
$l$ 
 is a {\bf linear} functional. 

 In addition I will take 
\begin{eqnarray*}
 M\frac{\partial^2 q}{\partial t^2}
 &=&
 F_{pf}(q,Q,t)   
\end{eqnarray*}
%

%%%%%\newpage 

 In the previous section 
$\delta_{\alpha,t,x_0}$
 was a function of 
$x$ 
 and it did not depend on 
$t$ 
$$
 \delta_{\alpha, t}(x) = \delta(x-x_0) \,;
$$
 As for
$l$,
 it was such that 
 for any function 
$F$  
 of 
$x$ 
 one had set 
$$
  <l|F> = F(x_0)
$$
 In addition we had there taken 
$$
 U_{0}(t,Q) = 4\gamma_c \Big(Q(t)-q(t)\Big) \,, \gamma_c=\gamma c \,,
$$
$$
 F_{pf}(q,Q,t)/M = -\Omega^2 \Big(q(t)-Q(t)\Big) + f_0(t) \,. 
$$

 Now we are mimicking, imitating arguments of the previous section:  
$$
 u(t) = u_{0}(t) 
        + \int_s^t V_{12}(t,\tau)f(\tau)d\tau 
$$

$$
 u(t) = u_{0}(t) 
        + \int_s^t V_{12}(t,\tau)
            \Big(-\delta_{\alpha,\tau,x_0}U_{0}(\tau,Q) +f_1(\tau)\Big)d\tau 
$$

$$
 u(t) = u_{0}(t) 
        + \int_s^t V_{12}(t,\tau)
            f_1(\tau)d\tau 
         - \int_s^t V_{12}(t,\tau)
            \Big(\delta_{\alpha,\tau,x_0}U_{0}(\tau,Q) \Big)d\tau 
$$

$$
 u_{01}(t) 
       := u_{0}(t) 
          + \int_s^t V_{12}(t,\tau)
              f_1(\tau)d\tau 
$$

$$
 u(t) = u_{01}(t) 
         - \int_s^t V_{12}(t,\tau)
            \Big(\delta_{\alpha,\tau,x_0}U_{0}(\tau,Q) \Big)d\tau 
$$

$$
 <l|u(t)> = <l|u_{01}(t)> 
           - <l|\int_s^t V_{12}(t,\tau)
            \Big(\delta_{\alpha,\tau,x_0}U_{0}(\tau,Q) \Big)d\tau >
$$

 We will suppose
$$
  <l|\int_s^t 
    V_{12}(t,\tau)\Big(\delta_{\alpha,\tau,x_0}U_{0}(\tau,Q) \Big)d\tau >
    = 
  \int_s^t <l|V_{12}(t,\tau)
            \Big(\delta_{\alpha,\tau,x_0}U_{0}(\tau,Q) \Big)>d\tau 
$$
 Hence 

$$
 <l|u(t)> = <l|u_{01}(t)> 
           - \int_s^t <l| V_{12}(t,\tau)
              \delta_{\alpha,\tau,x_0}>U_{0}(\tau,Q)d\tau 
$$
 Recall
$$
 Q(t) := <l|u(t)> 
$$
 and put 
$$
 K_{0}(t,\tau) := <l|V_{12}(t,\tau)\delta_{\alpha, \tau}> 
 \,, \qquad 
 Q_{0}(t) := <l|u_{0}(t)> 
 \,, \qquad 
$$
$$
 Q_{01}(t) := <l|u_{01}(t)>
              = <l|u_{0}(t)> 
                +\int_s^t <l|V_{12}(t,\tau)f_1(\tau)> d\tau 
$$
 Then obtain 
$$
 Q(t) =  Q_{0}(t) 
        -\int_s^t K_{0}(t,\tau)
                      U_{0}\Big(\tau, Q \Big) d\tau 
$$
 Let us now define an operator 
$\hat K_{0}$: 
$$
 ( \hat K_{0} h )(t) := \int_s^t K_{0}(t,\tau) h(\tau)d\tau 
$$
 Then the recent relation becomes as following: 
$$          
  Q =Q_{01} - \hat K_{0} \Big( U_{0}\Big(\cdot, Q \Big)\Big) 
$$
  If we take into account that 
$$
 M\frac{\partial^2 q}{\partial t^2}
 =
 F_{pf}(q,Q,t)   
$$
 then the equations of motion become  
\begin{eqnarray*}
 Q 
 &=&
 Q_{01} - \hat K_{0}\Big( U_{0}\Big(\cdot, Q \Big)\Big) 
 \\
 M\frac{\partial^2 q}{\partial t^2}
 &=&
 F_{pf}(q,Q,t)   
\end{eqnarray*}

 The character of this relation array is too general. 
 Little can be said about its properties without any specifying. 
 So, let us restrict ourselves: 
 First, we will assume 
\begin{eqnarray*}
 U_{0}(t,Q)    &=& 4\gamma_c \Big(Q(t)-q(t)\Big) \\ 
 F_{pf}(q,Q,t) &=& -K(q-Q)+F_0(t) = M\Big(-\Omega^2(q-Q)+f_0(t)\Big) 
\end{eqnarray*}
  and denote
$$
 \ddot q := \frac{\partial^2 q}{\partial t^2} 
$$
 Then the equations of motion become  
\begin{eqnarray*}
 Q 
 &=&
 Q_{01} - 4\gamma_c\hat K_{0} \Big( Q-q \Big)
 \\
 \ddot q 
 &=&
 -\Omega^2(q-Q)+f_0 
\end{eqnarray*}
 We infer: 

 \medskip

 1st step,
\begin{eqnarray*}
 Q + 4\gamma_c\hat K_{0} Q 
 &=&
 Q_{01} + 4\gamma_c\hat K_{0}q 
 \\
 \ddot q 
 &=&
 -\Omega^2(q-Q)+f_0 
\end{eqnarray*}

 2nd step,
\begin{eqnarray*}
 \Big(1 + 4\gamma_c\hat K_{0}\Big)Q 
 &=&
 Q_{01} + 4\gamma_c\hat K_{0}q 
 \\
 \ddot q 
 &=&
 -\Omega^2 q + \Omega^2 Q + f_0 
\end{eqnarray*}

 3rd step, 
\begin{eqnarray*}
 \Big(1 + 4\gamma_c\hat K_{0}\Big)Q 
 &=&
 Q_{01} + 4\gamma_c\hat K_{0}q 
 \\
% \qquad 
 \Big(1 + 4\gamma_c\hat K_{0}\Big) \ddot q 
 &=&
 \Big(1 + 4\gamma_c\hat K_{0}\Big)
 \Big( -\Omega^2 q + \Omega^2 Q + f_0 \Big)  
\end{eqnarray*}

 4th step, 
\begin{eqnarray*}
 \Big(1 + 4\gamma_c\hat K_{0}\Big)Q 
 &=&
 Q_{01} + 4\gamma_c\hat K_{0}q 
 \\
 \Big(1 + 4\gamma_c\hat K_{0}\Big) \ddot q 
 &=& 
 -\Big(1 + 4\gamma_c\hat K_{0}\Big)\Omega^2 q 
 +\Big(1 + 4\gamma_c\hat K_{0}\Big)\Omega^2 Q 
 +\Big(1 + 4\gamma_c\hat K_{0}\Big)f_0 
 \\ 
 &=& 
 -\Omega^2 \Big(1 + 4\gamma_c\hat K_{0}\Big)q 
 +\Omega^2 \Big(1 + 4\gamma_c\hat K_{0}\Big)Q 
 +\Big(1 + 4\gamma_c\hat K_{0}\Big)f_0 
 \\ 
 &&\makebox[41ex][l]{ \hfill$
        \left\{
        \begin{array}{r}\mbox{ replace here }
                        \Big(1 + 4\gamma_c\hat K_{0}\Big)Q 
                      \\
                         \mbox{ by }
                         \Big( Q_{01} + 4\gamma_c\hat K_{0}q \Big) 
         \end{array}\right\}$}
 \\ 
 &=& 
 -\Omega^2 \Big(1 + 4\gamma_c\hat K_{0}\Big)q 
 +\Omega^2 \Big( Q_{01} + 4\gamma_c\hat K_{0}q \Big) 
 +\Big(1 + 4\gamma_c\hat K_{0}\Big)f_0 
 \\ 
 &=& 
 -\Omega^2 q 
 +\Omega^2 Q_{01} 
 +\Big(1 + 4\gamma_c\hat K_{0}\Big)f_0 
 \\ 
\end{eqnarray*}
 We have now seen: 
$$
 \Big(1 + 4\gamma_c\hat K_{0}\Big) \ddot q 
 =
 -\Omega^2 q 
 +\Omega^2 Q_{01} 
 +\Big(1 + 4\gamma_c\hat K_{0}\Big)f_0 
$$
 It means 
\begin{eqnarray*}
\makebox[5ex][l]{$\displaystyle 
 \frac{\partial^2 q(\tau)}{\partial \tau^2}
 +4\gamma_c 
 \int_s^t K_0(t,\tau) \frac{\partial^2 q(\tau)}{\partial \tau^2}d\tau  
  $}%% end of the makebox 
\\&=& 
  -\Omega^2 q(t) + \Omega^2 Q_{01}(t)  
%\\&&{}
 + f_0(t)
 +4\gamma_c \int_s^t K_0(t,\tau)
 \Big( f_0(\tau)\Big) d\tau 
\end{eqnarray*}
 We rewrite this equation as following: 
$$
  \frac{\partial^2 q}{\partial t^2}
 =
 -4\gamma_c 
 \int_s^t K_0(t,\tau) \frac{\partial^2 q(\tau)}{\partial \tau^2}d\tau  
 -\Omega^2 q 
 +\Omega^2 Q_{01}(t)
 +f_0(t) +4\gamma_c \int_s^t K_0(t,\tau)f_0(\tau) d\tau 
$$ 
 and then, replacing 
$Q_{01}(t)$
 by 
$<l|u_{01}(t)>$ ,  
 as 
\begin{eqnarray*} 
\makebox[28ex][l]{$\displaystyle 
  \frac{\partial^2 q}{\partial t^2} 
 =
 -4\gamma_c 
 \int_s^t K_0(t,\tau) \frac{\partial^2 q(\tau)}{\partial \tau^2}d\tau  
 -\Omega^2 q 
  $}%% end of the makebox 
\\[-1\bigskipamount]
 &&{ } +\Omega^2 <l|u_{01}(t)> 
 +f_0(t) +4\gamma_c \int_s^t K_0(t,\tau)f_0(\tau) d\tau 
\end{eqnarray*}

 Thus, we have obtained an {\bf insulated } equation for 
$q$, 
 and the equation resembles ones in the previous section. 
 To continue, we must specify the nature of the abstract terms, of 
$u(t)$. 

 If they, 
$u(t)$, 
 all are elements of an functional space, 
 say simple, 
 if 
$u(t)$ 
 are functions of a {`spacial'} variable 
$x$, 
 our constructions and formulae 
 appear 
 as follows:

%%%%%%%%%%%%\newpage 

$$
 u(t,x) = u_{0}(t,s,x)
        + \int_s^t \Big(V_{12}(t,\tau)f(\tau)\Big)(x)d\tau 
$$

$$
 u(t,x) = u_{0}(t,s,x)
        + \int_s^t \int_X V_{12}(t,\tau,x,\xi)f(\tau,\xi)d\xi d\tau 
$$
 In this case we have 
$$
f(t,x) :=-\delta_{\alpha,t,x_0}(x)U_{0}(t,Q) + f_1(t,x) 
$$
 ( in the previous section, 
$ \delta_{\alpha,t,x_0}(x,x_0) = \delta(x-x_0) $ 
 )
 Then 
$$
 u(t,x) = u_{0}(t,s,x)
         -\int_s^t \int_X V_{12}(t,\tau,x,\xi) 
             \Big(\delta_{\alpha,\tau,x_0}(\xi)U_{0}(\tau,Q) 
                   - f_1(\tau,\xi)\Big) d\xi d\tau 
$$
 Hence
\begin{eqnarray*}
 <l|u(t)>
 &=& <l|u(t,\cdot)> 
\\
 &=& <l|u_{0}(t,s,\cdot)
            -\int_s^t \int_X V_{12}(t,\tau,\cdot,\xi) 
             \Big(\delta_{\alpha,\tau,x_0}(\xi)U_{0}(\tau,Q) 
                   - f_1(\tau,\xi)\Big) d\xi d\tau >
\\
 &=& <l|u_{0}(t,s,\cdot)>
            -<l|\int_s^t \int_X V_{12}(t,\tau,\cdot,\xi) 
              \delta_{\alpha,\tau,x_0}(\xi)
              U_{0}(\tau,Q)d\xi d\tau> 
\\&&{}           +<l|\int_s^t \int_X V_{12}(t,\tau,\cdot,\xi)
                    f_1(\tau,\xi) d\xi d\tau> 
\\
 &=& <l|u_{0}(t,s,\cdot)>
            -\int_s^t <l|\int_X V_{12}(t,\tau,\cdot,\xi) 
              \delta_{\alpha,\tau,x_0}(\xi)d\xi>
              U_{0}(\tau,Q) d\tau 
\\&&{}           +<l|\int_s^t \int_X V_{12}(t,\tau,\cdot,\xi)
                    f_1(\tau,\xi) d\xi d\tau> 
\end{eqnarray*}
 Hence
\begin{eqnarray*}
 <l|u(t)>
 &=& <l|u(t,\cdot)> 
\\
 &=& <l|u_{0}(t,s,\cdot)>
            -\int_s^t <l|\int_X V_{12}(t,\tau,\cdot,\xi) 
              \delta_{\alpha,\tau,x_0}(\xi)d\xi>
              U_{0}(\tau,Q) d\tau 
\end{eqnarray*}
 i.e. 
\begin{eqnarray*}
 Q(t) 
 &=& Q_{01}(t)
            -\int_s^t <l|\int_X V_{12}(t,\tau,\cdot,\xi) 
              \delta_{\alpha,\tau,x_0}(\xi)d\xi>
              U_{0}(\tau,Q) d\tau 
\end{eqnarray*}
\begin{eqnarray*}
 Q(t) 
 &=& Q_{01}(t)
            -\int_s^t K_0(t,\tau) 
              U_{0}(\tau,Q) d\tau 
\end{eqnarray*}
 Of course 
$$
   K_0(t,\tau) 
  := <l|\int_X V_{12}(t,\tau,\cdot,\xi)\delta_{\alpha,\tau,x_0}(\xi) d\xi > 
$$
$$
 Q_{01} =  <l|u_{0}(t,s,\cdot)>
                +<l|\int_s^t \int_X V_{12}(t,\tau,\cdot,\xi)
                    f_1(\tau,\xi) d\xi d\tau> 
$$

%%%%%%%%%%%%%\newpage 

 If we presume  
$B,l$ 
 to be the same as before, in the previous section, 
 i.e., if  
$B$ 
 is  
$c^2\frac{\partial^2}{\partial x^2}$
 on the whole line, 
 and 
$ <l|F> = F(x_0) $, 
then 

\begin{eqnarray*} 
\makebox[5ex][l]{$\displaystyle 
 \int_s^t \int_X V_{12}(t,\tau,x,\xi) f(\tau ,\xi) d\xi d\tau 
  $}%% end of the makebox 
\\&=& 
 \frac12
 \int_s^t 
 \Big( \tilde f(\tau,x+c(t-\tau))-\tilde f(\tau,x-c(t-\tau))\Big) d\tau 
\end{eqnarray*} 
\begin{eqnarray*} 
\makebox[5ex][l]{$\displaystyle 
 \int_s^t \int_X V_{12}(t,\tau,x,\xi) 
 \delta_{\alpha,\tau,x_0}(\xi)U_{0}(\tau,Q) 
 d\xi d\tau 
  $}%% end of the makebox 
\\&=& 
 \frac{1}{2c}
 \int_s^t 
 \Big( 
 \tilde \delta_{\alpha,\tau,x_0}(x+c(t-\tau))
  -\tilde\delta_{\alpha,\tau,x_0}(x-c(t-\tau))
 \Big) 
 U_{0}(\tau,Q) d\tau 
\end{eqnarray*} 
 Thus we obtain 
\begin{eqnarray*} 
 K_0(t,\tau) 
 &=& 
 \frac{1}{2c}
 \Big( 
 \tilde \delta_{\alpha,\tau,x_0}(x_0+c(t-\tau))
  -\tilde\delta_{\alpha,\tau,x_0}(x_0-c(t-\tau))
 \Big) 
\end{eqnarray*} 
 Here  
$\tilde \delta_{\alpha,\tau,x_0}$ 
 stands for any primitive of 
$\delta_{\alpha,\tau,x_0}$, 
 i.e.,
$$
  \frac{\partial\tilde \delta_{\alpha,\tau,x_0}(x)}{\partial x} 
  = \delta_{\alpha,\tau,x_0}(x)
$$
 Note once again that
$$
 \tilde \delta_{\alpha,\tau,x_0}(x+c(t-\tau)) 
 -\tilde \delta_{\alpha,\tau,x_0}(x-c(t-\tau)) 
$$
 does not depend on whatever primitive of 
$\delta_{\alpha,\tau,x_0}$ 
 which one has chosen!!!  

 Suppose, 
$\delta_{\alpha,t,x_0}(x)$ 
 is of the form 
$$
 \delta_{\alpha,t,x_0}(x) = \delta_{\alpha}(x-x_0)
$$ 
 Then we obtain 
\begin{eqnarray*} 
 K_0(t,\tau) 
 &=& 
 \frac{1}{2c}
 \Big( 
 \tilde \delta_{\alpha,\tau,x_0}(x_0+c(t-\tau))
  -\tilde\delta_{\alpha,\tau,x_0}(x_0-c(t-\tau))
 \Big) 
\\&=& 
 \frac{1}{2c}
 \Big( 
 \tilde \delta_{\alpha}(c(t-\tau))-\tilde \delta_{\alpha}(c(\tau-t))
 \Big) 
\end{eqnarray*} 
  
 In this case, we have finally seen,
\begin{eqnarray*}
  \frac{\partial^2 q(t)}{\partial t^2}
 &=&
 -2\frac{\gamma_c}{c} 
 \int_s^t 
 \Bigl( 
 \tilde \delta_{\alpha}(c(t-\tau))-\tilde \delta_{\alpha}(c(\tau-t))
 \Bigr) 
                   \frac{\partial^2 q(\tau)}{\partial \tau^2}d\tau  
 -\Omega^2 q 
 +\Omega^2 Q_{01}(t)
 \\&&{ }\qquad\qquad\qquad 
 +f_0(t) 
 +2\frac{\gamma_c}{c} 
 \int_s^t  
 \Bigl( 
 \tilde \delta_{\alpha}(c(t-\tau))-\tilde \delta_{\alpha}(c(\tau-t))
 \Bigr) 
 f_0(\tau) d\tau 
 \\&=&
 -2\gamma 
 \int_s^t 
 \Bigl( 
 \tilde \delta_{\alpha}(c(t-\tau))-\tilde \delta_{\alpha}(c(\tau-t))
 \Bigr) 
                   \frac{\partial^2 q(\tau)}{\partial \tau^2}d\tau  
 -\Omega^2 q 
 +\Omega^2 Q_{01}(t)
 \\&&{ }\qquad\qquad\qquad 
 +f_0(t) 
 +2\gamma 
 \int_s^t  
 \Bigl( 
 \tilde \delta_{\alpha}(c(t-\tau))-\tilde \delta_{\alpha}(c(\tau-t))
 \Bigr) 
 f_0(\tau) d\tau 
\end{eqnarray*}
 In any case, 
$$
  \frac{\partial^2 q(t)}{\partial t^2}-f_0(t) 
  +2\frac{\gamma_c}{c} \int_s^t 2K_0(t,\tau) 
 \Bigl( 
 \frac{\partial^2 q(\tau)}{\partial \tau^2}-f_0(\tau)
 \Bigl) 
 d\tau  
 = 
 -\Omega^2\Bigl(q(t) - Q_{01}(t)\Bigr) 
$$ 
 or, more generally, 
\begin{eqnarray*}
 Q(t) 
 &=& Q_{01}(t)
            -\int_s^t K_0(t,\tau)  
              U_{0}(\tau,Q) d\tau \,,
 \\
 M\frac{\partial^2 q(t)}{\partial t^2}
 &=&
 F_{pf}(q,Q,t)  \,.  
\end{eqnarray*}

%%%%%%%%%%%
%%%%%%%%%%%

\newpage 
\section{APPENDIX}
\subsection{ APPENDIX A. Abstract Linear Response Formula } 

 We recall some abstract linear response formulae.

 First, we recall linear response formula for 
 the {\bf first} order differential equation. 
 We indicate only formula itself, minimum of details. 
 We say only, a deriving of this formula is based on an abstract modification 
 of the method of {\bf variation of constant}. 
 As for detailed description, how one produces such a formula, 
 as well as for proofs, background and all that,
 see suitable standard manuals.

 Let 
$\{A(t)\}_t$ 
 stand for a family of linear operators, 
 so that the {\bf homogeneous } linear equation
$$
 \frac{\partial v(t)}{\partial t}
 =A(t)v(t),\quad  v(t)\Big|_{t=s} = v_{initial}, \qquad    (t\geq s)
$$
 has a unique solution 
$v(t)$, whatever `initial time' $s$ and initial data $v_{initial}$ may be.

 Let
$$
 V_{t,s}
$$
 denote the propagator, alias evolution operator, to the former system, i.e.,
$$
 v(t) =  V_{t,s}v_{initial} =  V_{t,s}v(s) .
$$ 
 One can show that 
$$
  V_{t,r}V_{r,s} = V_{t,s} \,, \mbox{ if } t\geq r\geq s,  
  \makebox[0pt][l]{ ( consistency relation )}
$$
$$
  V_{t,t} = V_{r,r} = V_{s,s} = I \,. 
$$
 The linear response formula reads:
 The formula
$$
 \fbox{$\displaystyle
  u(t) = V_{t,s}u_{initial} + \int_s^t V_{t,\tau}f(\tau)d\tau $}
$$ 
 gives a solution to the {\bf inhomogeneous} equation 

$$
 \fbox{$\displaystyle
 \frac{\partial u(t)}{\partial t}
 =A(t)u + f(t), \quad u(t)\Big|_{t=s} \equiv u(s) = u_{initial} .$}
$$
 We will write the response formula as following:
$$
 \fbox{$\displaystyle
  u(t) = V_{t,s}u(s) + \int_s^t V_{t,\tau}f(\tau)d\tau .$}
$$ 
 One can now produce the response formula 
 for the {\bf second} order differential equations. 
 As befor, the response formula accommodates the formula,
 which describes solution to a {\bf homogeneous} equation, to the case,
 where one need to solve the associated {\bf inhomogeneous} one.

\newpage 
 Now then. When solving the equation
$$
  \frac{\partial^2 u(t)}{\partial t^2} 
 = K(t)\frac{\partial u(t)}{\partial t} + B(t)u(t)+f(t) 
$$
 denote
$$
 u_1(t):=\frac{\partial u(t)}{\partial t} \,.
$$
 Write the equation to be solved as
\begin{eqnarray*}
 \frac{\partial }{\partial t} 
\Big( \begin{array}{cc}
 u(t) \\ u_1(t)
\end{array}\Big)
 &=& 
\Big( \begin{array}{cc}
 u_1(t) \\ K(t)u_(t) +B(t)u(t) + f(t)
\end{array}\Big)
 \\&&{ } 
 = 
\Big( \begin{array}{cc}
 0 & 1 \\ B(t) & K(t) \\ 
\end{array}\Big)
\Big( \begin{array}{cc}
 u(t) \\ u_1(t) 
\end{array}\Big)
 + 
\Big( \begin{array}{cc}
 0 \\ f(t) 
\end{array}\Big)
\end{eqnarray*}
 Let
$$
 V_{t,s}
$$
 denote the propagator to 
$$
 \frac{\partial }{\partial t} 
\Big( \begin{array}{cc}
 v(t) \\ v_1(t) 
\end{array}\Big)
 = 
\Big( \begin{array}{cc}
 0 & 1 \\ B(t) & K(t) \\
\end{array}\Big)
\Big( \begin{array}{cc}
 v(t) \\ v_1(t) 
\end{array}\Big)
$$
 Then
$$
\Big( \begin{array}{cc}
 u(t) \\ u_1(t) 
\end{array}\Big)
 = V_{t,s}
\Big( \begin{array}{cc}
 u(s) \\ u_1(s) 
\end{array}\Big)
 + \int_s^t V_{t,\tau} 
\Big( \begin{array}{cc}
 0 \\ f(\tau) 
\end{array}\Big)
 d\tau 
$$
 Apply usual matrix-form representation
$$
 V_{t,s} = 
\Big( \begin{array}{cc}
 V_{11}(t,s) & V_{12}(t,s) \\  V_{21}(t,s) & V_{22}(t,s) 
\end{array}\Big)
$$
 Then obtain finally 

$$
\fbox{
\parbox{\textwidth}
{
\begin{eqnarray*} 
 u(t) &=& V_{11}(t,s)u(s) 
        + V_{12}(t,s)\Big(\frac{\partial u(t)}{\partial t} \Big|_{t=s} \Big) 
        + \int_s^t V_{12}(t,\tau)f(\tau)d\tau 
\end{eqnarray*}
}%%end of parbox
} %% end of fbox
$$

 Notice,
$$
\Big( \begin{array}{cc}
 V_{11}(t,s) & V_{12}(t,s) \\  V_{21}(t,s) & V_{22}(t,s) 
\end{array}\Big)
 =
\Big( \begin{array}{cc}
 I & 0 \\ 0  & I 
\end{array}\Big)
$$
 and the formula 

$$
 v(t) = V_{11}(t,s)v(s) 
        + V_{12}(t,s)\Big(\frac{\partial v(t)}{\partial t} \Big|_{t=s} \Big) 
$$
 is a formula which describes solutions to the homogeneous equation. 
 Sometimes, it is instructive to notice that 
$$
\frac{\partial }{\partial t} v(t)
 = \frac{\partial }{\partial t} V_{11}(t,s)v(s) 
        + \frac{\partial }{\partial t} V_{12}(t,s)
             \Big(\frac{\partial v(t)}{\partial t} \Big|_{t=s} \Big) 
$$
 and hence 
$$
 V_{21}(t,s) = \frac{\partial }{\partial t} V_{11}(t,s)
 \,, \quad 
 V_{22}(t,s) = \frac{\partial }{\partial t} V_{12}(t,s)
 \,. 
$$

\newpage 
\subsection{ APPENDIX B.}
\subsection*%
{ Response Formula for $\ddot y(t) +2\gamma\dot y(t) + \Omega^2 y(t) = f(t)$} 

 If one applies the response formula to the equation in the title, 
 then one finds: 
\begin{eqnarray*} 
 y(t) 
 &=& e^{-\gamma (t-s)}\Bigl(\cos\Omega_\gamma (t-s)
   +\gamma\frac{\sin\Omega_\gamma (t-s)}{\Omega_\gamma }\Bigr) y(s) 
 \\&&{ } \hspace*{5ex}
  + e^{-\gamma (t-s)}\frac{\sin\Omega_\gamma (t-s)}{\Omega_\gamma } 
                \Big(\frac{\partial y(t)}{\partial t} \Big|_{t=s} \Big) 
 \\&&{ } \hspace*{10ex}
  + \int_s^t 
           e^{-\gamma (t-\tau)}
            \frac{\sin\Omega_\gamma (t-\tau)}{\Omega_\gamma }f(\tau)d\tau 
\end{eqnarray*}
 where 
$\Omega_\gamma ^2 = \Omega^2 -\gamma ^2$
\footnote{
 if 
$\Omega_\gamma =0$,
 then read 
$\frac{\sin\Omega_\gamma t}{\Omega_\gamma } = t$, 
 of course}.

 Actually, 
 the solution to the joined homogeneous equation
$$
\ddot x(t) +2\gamma \dot x(t) + \Omega^2 x(t) = 0 
$$
 is this:
\begin{eqnarray*}
 x(t) 
 &=& e^{-\gamma t}
      \Bigl(\cos\Omega_\gamma t
            +\gamma\frac{\sin\Omega_\gamma t }{\Omega_\gamma }
            \Bigr)x(0) 
     + e^{-\gamma t}\frac{\sin\Omega_\gamma t}{\Omega_\gamma } \dot x(0)
 \\ 
 \dot 
 x(t) 
 &=& e^{-\gamma t}
       \Bigl(-{\Omega_\gamma}^2 -\gamma^2\Bigr)
       \Bigl(\frac{\sin\Omega_\gamma t}{\Omega_\gamma }\Bigr) x(0) 
     +e^{-\gamma t}
       \Bigl(\cos\Omega_\gamma t
         -\gamma\frac{\sin\Omega_\gamma t }{\Omega_\gamma }\Bigr)\dot x(0) 
\end{eqnarray*}
  Since 
$\Omega_\gamma ^2 = \Omega^2 -\gamma ^2$, 
  we can besides write  
\begin{eqnarray*}
 x(t) 
 &=& e^{-\gamma t}
      \Bigl(\cos\Omega_\gamma t
            +\gamma\frac{\sin\Omega_\gamma t }{\Omega_\gamma }
            \Bigr)x(0) 
     + e^{-\gamma t}\frac{\sin\Omega_\gamma t}{\Omega_\gamma } \dot x(0)
 \\ 
 \dot 
 x(t) 
 &=& e^{-\gamma t}
       \Bigl( -{\Omega}^2 \Bigr)
       \Bigl(\frac{\sin\Omega_\gamma t}{\Omega_\gamma }\Bigr) x(0) 
     +e^{-\gamma t}
       \Bigl(\cos\Omega_\gamma t
         -\gamma\frac{\sin\Omega_\gamma t }{\Omega_\gamma }\Bigr)\dot x(0) 
\end{eqnarray*}
 The coefficients of the latter equation for 
$x(t)$  
 are constant. Hence  
\begin{eqnarray*}
\makebox[30ex][l]{$\displaystyle 
 x(t) 
 = e^{-\gamma (t-s)}
 \Bigl(\cos\Omega_\gamma (t-s)
     +\gamma\frac{\sin\Omega_\gamma (t-s)}{\Omega_\gamma }\Bigr)x(s) 
  $}%end of makebox
 \\&& + e^{-\gamma (t-s)}
   \frac{\sin\Omega_\gamma (t-s)}{\Omega_\gamma } \dot x(s)
 \\ 
\makebox[30ex][l]{$\displaystyle 
 \dot 
 x(t) 
 = e^{-\gamma (t-s)}
       \Bigl( -{\Omega}^2 \Bigr)
         \frac{\sin\Omega_\gamma (t-s)}{\Omega_\gamma }
     x(s) 
  $}%end of makebox
  \\&& +e^{-\gamma (t-s)}
       \Big(\cos\Omega_\gamma (t-s) 
          -\gamma\frac{\sin\Omega_\gamma (t-s)}{\Omega_\gamma }
                                \Bigr)\dot x(s) 
\end{eqnarray*}
 Thus we have seen, that the formula for the propagator 
$V_{t,s}$, 
 which associates with the homogeneous equation 
$\ddot x(t) +2\gamma \dot x(t) + \Omega^2 x(t) = 0$, 
 is this:
\begin{eqnarray*}
 V_{t,s} &=& 
\left( \begin{array}{cc}
 V_{11}(t,s) & V_{12}(t,s) \\  V_{21}(t,s) & V_{22}(t,s) 
\end{array}\right)
 \\&&{ }
 = 
\left(
\begin{array}{cc}
  e^{-\gamma (t-s)}
  \Bigl(\cos\Omega_\gamma (t-s)
        +\gamma\frac{\sin\Omega_\gamma (t-s)}{\Omega_\gamma }
                                \Bigr) 
 & e^{-\gamma (t-s)}\frac{\sin\Omega_\gamma (t-s)}{\Omega_\gamma } 
 \\ e^{-\gamma (t-s)}
       \Bigl( -{\Omega}^2 \Bigr)
       \frac{\sin\Omega_\gamma (t-s)}{\Omega_\gamma }
         
 & e^{-\gamma (t-s)}
  \Big(\cos\Omega_\gamma (t-s)
        -\gamma\frac{\sin\Omega_\gamma (t-s)}{\Omega_\gamma }\Bigr) 
\end{array}
\right) 
\end{eqnarray*}

 Repeat, finally, the abstract response formula 
$$
\fbox{
\parbox{\textwidth}
{
\begin{eqnarray*} 
 u(t) &=& V_{11}(t,s)u(s) 
        + V_{12}(t,s)\Big(\frac{\partial u(t)}{\partial t} \Big|_{t=s} \Big) 
        + \int_s^t V_{12}(t,\tau)f(\tau)d\tau 
\end{eqnarray*}
}%%end of parbox
} %% end of fbox
$$
 and replace there the proper terms according to 
$$
 V_{11}(t,s) 
 = e^{-\gamma (t-s)}\Bigl(\cos\Omega_\gamma (t-s)
   +\gamma\frac{\sin\Omega_\gamma (t-s)}{\Omega_\gamma }\Bigr) 
$$
$$
 V_{12}(t,s) 
 = e^{-\gamma (t-s)}\frac{\sin\Omega_\gamma (t-s)}{\Omega_\gamma } 
$$

\newpage 
\subsection{ APPENDIX C.}
\subsection*%
 { Particular Solution to}
\subsection*%
 { $\ddot y(t) +2\gamma\dot y(t) + \Omega^2 y(t) 
   = A_s\sin(kt+\phi)+A_c\cos(kt+\phi) $} 

 We have 

\begin{eqnarray*} 
 \makebox[21ex][l]{$\displaystyle 
 \Bigl(\frac{\partial^2 }{\partial t^2} 
 +2\gamma\frac{\partial }{\partial t} 
 +\Omega^2\Bigr)\Bigl(a_s\sin(kt+\phi)+a_c\cos(kt+\phi)\Bigr) 
$}% end makebox
 \\ &=& 
 \Bigl(-k^2+\Omega^2\Bigl)a_s\sin(kt+\phi) 
     +2\gamma k a_s \cos(kt+\phi) 
 \\ &&{ } 
     -2\gamma k a_c\sin(kt+\phi) 
         + \Bigl(-k^2+\Omega^2\Bigl) a_c \cos(kt+\phi) 
\end{eqnarray*}
  In order to
$$
 \left( 
 \begin{array}{cc} 
   -k^2+\Omega^2  & -2\gamma k 
 \\
   2\gamma k      & -k^2+\Omega^2 
 \end{array} 
 \right) 
\left( 
 \begin{array}{cc} 
  a_s  \\ a_c 
 \end{array} 
 \right) 
 = 
 \left( 
 \begin{array}{cc} 
  A_s  \\ A_c 
 \end{array} 
 \right) 
$$
 it is sufficient that
$$
 \Bigl((-k^2+\Omega^2)^2 + (2\gamma k)^2 \Bigl) 
 \left( 
 \begin{array}{cc} 
  a_s  \\ a_c 
 \end{array} 
 \right) 
 = 
 \left( 
 \begin{array}{cc}
   -k^2+\Omega^2  & 2\gamma k 
 \\ 
   -2\gamma k     & -k^2+\Omega^2 
 \end{array} 
 \right) 
 \left( 
 \begin{array}{cc} 
  A_s  \\ A_c 
 \end{array} 
 \right) 
$$
 Thus, if 
$$
 (-k^2+\Omega^2)^2 + (2\gamma k)^2 \not= 0 
$$
 then
$$
 \left( 
 \begin{array}{cc} 
  a_s  \\ a_c 
 \end{array} 
 \right) 
 = 
 \frac{1}{ (-k^2+\Omega^2)^2 + (2\gamma k)^2 } 
 \left( 
 \begin{array}{cc}
   -k^2+\Omega^2  & 2\gamma k 
 \\ 
  - 2\gamma k     & -k^2+\Omega^2 
 \end{array} 
 \right) 
 \left( 
 \begin{array}{cc} 
  A_s  \\ A_c 
 \end{array} 
 \right) 
$$

\newpage 
\subsection{ APPENDIX D. Helpful formulae } 
%\indent
 Derivatives. 
 1) 

\begin{eqnarray*}
\makebox[8ex][l]{$\displaystyle 
 \frac{\partial 
 e^{-\gamma t}
 \Bigl(
 A\cos(\omega t+\varphi)+B\sin(\omega t +\varphi)
 \Bigr)}{\partial t}
$}%%end of makebox
\\ &=& -\gamma e^{-\gamma t}
 \Bigl(
 A\cos(\omega t+\varphi)+B\sin(\omega t+\varphi)\Bigr) 
\\ &&{ }
 +e^{-\gamma t}
 \Bigl(-\omega A\sin(\omega t+\varphi)+\omega B\cos(\omega t+\varphi)\Bigr)
\\ &=&
 e^{-\gamma t}\Bigl(
 (-\gamma A+\omega B)\cos(\omega t+\varphi)
 +(-\omega A-\gamma B)\sin(\omega t+\varphi)
  \Bigr) 
\end{eqnarray*}
 2)
\begin{eqnarray*}
 \frac{\partial 
 e^{-\gamma t}
 \Bigl(
 -\gamma \cos(\omega t+\varphi)+\omega\sin(\omega t +\varphi)
 \Bigr)}{\partial t}
 =
 e^{-\gamma t}(\gamma ^2+\omega^2)\cos(\omega t+\varphi)
\end{eqnarray*}
 3)
\begin{eqnarray*}
 \frac{\partial 
 e^{-\gamma t}
 \Bigl(
 \omega\cos(\omega t+\varphi)+\gamma \sin(\omega t +\varphi)
 \Bigr)}{\partial t}
 =
 -e^{-\gamma t}(\gamma ^2+\omega^2)\sin(\omega t+\varphi)
\end{eqnarray*}
 2a)
\begin{eqnarray*}
\makebox[30ex][l]{$\displaystyle 
 \frac{\partial 
 e^{-\gamma (t-\tau)}
 \Bigl(
 -\gamma \cos(\omega(t-\tau)+\varphi)+\omega\sin(\omega(t-\tau) +\varphi)
 \Bigr)}{\partial \tau}
$}%%end of makebox
 \\&&{ } 
 =
 -e^{-\gamma (t-\tau)}(\gamma ^2+\omega^2)\cos(\omega(t-\tau)+\varphi)
\end{eqnarray*}
 3a)
\begin{eqnarray*}
\makebox[30ex][l]{$\displaystyle 
 \frac{\partial 
 e^{-\gamma t}
 \Bigl(
 \omega\cos(\omega(t-\tau)+\varphi)+\gamma \sin(\omega(t-\tau) +\varphi)
 \Bigr)}{\partial \tau}
$}%%end of makebox
 \\&&{ } 
 =
 e^{-\gamma (t-\tau)}(\gamma ^2+\omega^2)\sin(\omega(t-\tau)+\varphi)
\end{eqnarray*}
 Other.

$$
 2\sin\alpha\sin\beta = cos(\alpha-\beta) - cos(\alpha+\beta)
$$

\begin{eqnarray*}
 \omega(t-\tau)-\omega_s(\tau-t_0)
 &=&
 (\omega+\omega_s)(t-\tau)-\omega_s(t-\tau)-\omega_s(\tau-t_0)
\\&=&
 (\omega+\omega_s)(t-\tau)-\omega_s(t-t_0)
\end{eqnarray*}

\begin{eqnarray*}
 \omega(t-\tau)+\omega_s(\tau-t_0)
 &=&
 (\omega-\omega_s)(t-\tau)+\omega_s(t-\tau)+\omega_s(\tau-t_0)
\\&=&
 (\omega-\omega_s)(t-\tau)+\omega_s(t-t_0)
\end{eqnarray*}

\begin{eqnarray*}
\makebox[0ex][l]{$\displaystyle 
 2\sin\omega(t-\tau)\sin\omega_s(\tau-t_0)
$}%%end of makebox
 \\&&{ } 
 =
 \cos\Bigl((\omega+\omega_s)(t-\tau)-\omega_s(t-t_0)\Bigr)
 -\cos\Bigl((\omega-\omega_s)(t-\tau)+\omega_s(t-t_0)\Bigr)
\end{eqnarray*}

\newpage 

\bibliographystyle{unsrt}

\end{document}